\documentclass[a4paper]{article}
\usepackage{amssymb,amsmath,amscd,amsfonts,amsthm,latexsym,bbm}
\usepackage{xcolor,graphpap,fancybox,framed,enumerate,xcolor,mathtools}
\usepackage[english,activeacute]{babel}
\usepackage{needspace}
\usepackage{authblk}
\usepackage%[dvipdfmx]
{graphicx} 

\usepackage[right=2.5cm,left=2.5cm,top=2.5cm,bottom=2.5cm]{geometry}

\newtheorem{theorem}{Theorem}[section]
\newtheorem{prop}[theorem]{Proposition}
\newtheorem{coro}[theorem]{Corollary}

\newtheorem{remark}[theorem]{Remark}
\newenvironment{demo}{ \noindent \emph{\textbf{Proof:}}}{\hfill$\square$\\}

\newcommand{\RR}{\mathbb{R}}
\newcommand{\NN}{\mathbb{N}}
\newcommand{\CC}{\mathbb{C}}

\newcommand{\Bc}{\mathcal{B}}
\newcommand{\Cc}{\mathcal{C}}

\newcommand{\Hc}{\mathcal{H}}

\newcommand{\Lc}{\mathcal{L}}
\newcommand{\Mc}{\mathcal{M}}

\newcommand{\Xc}{\mathcal{X}}

\newcommand{\ra}{\rangle}
\newcommand{\la}{\langle}

\newcommand{\ddd}{\partial}

\newcommand{\iii}{{\, \vert\kern-0.25ex\vert\kern-0.25ex\vert\, }}
\newcommand{\grad}{\nabla}

\newcommand{\Tr}{\operatorname{Tr}}
\renewcommand{\div}{\operatorname{div}}
\newcommand{\Un}{\mathbbm{1}}  %{1\hspace{-1.5mm}1}
\newcommand{\no}{n$^{\text{o}}$}
\newcommand{\dd}{\,{\text{\rm d}}}
\newcommand{\Diff}{{\text{\rm Diff}}}
\newcommand{\meas}{{\text{\rm meas}}}
\newcommand{\id}{{\text{\rm id}}}

\newcommand{\pc}{ \usefont{T1}{cmtl}{m}{n} \selectfont}

\newdimen\texpscorrection
\newdimen\figcenter
\def\figurewithtex #1 #2 #3 #4 #5\cr{\null
  {\goodbreak\figcenter=\hsize\relax
  \advance\figcenter by -#4truecm
  \divide\figcenter by 2
  \begin{figure}[hbt]
  \vskip #3truecm\noindent\hskip\figcenter
  \includegraphics{#1}{\hskip\texpscorrection\input #2 }
  \vskip 0.8truecm{\baselineskip=0.8\baselineskip
  \noindent \vbox{\noindent {\footnotesize #5}}\par}
  \end{figure}}}
\def\point#1 #2 #3 {\rlap{\kern #1 truecm
\raise #2 truecm \hbox{#3}}}

\numberwithin{equation}{section}

\begin{document}

\title{{\bf Schr\"odinger equation in moving domains}}

\author[1]{Alessandro \textsc{Duca}}
\author[2]{Romain \textsc{Joly}}
\affil[1]{{\small Universit\'e Grenoble Alpes, CNRS, Institut Fourier, F-38000 Grenoble, France
{\pc alessandro.duca@univ-grenoble-alpes.fr}}}
\affil[2]{\small Universit\'e Grenoble Alpes, CNRS, Institut Fourier, F-38000 Grenoble, France {\pc romain.joly@univ-grenoble-alpes.fr}}

\date{}

\maketitle

\begin{abstract} 
We consider the Schr\"odinger equation
\begin{equation}\label{eq_abstract} 
i\partial_t u(t)=-\Delta u(t)~~~~~\text{ on }\Omega(t) \tag{$\ast$}
\end{equation}
where $\Omega(t)\subset\RR^N$ is a moving domain depending on the time $t\in [0,T]$. The aim
of this work is to provide a meaning to the solutions of such an equation. We 
use the existence of a bounded reference domain $\Omega_0$ and a specific family of unitary maps $h^\sharp(t): 
L^2(\Omega(t),\CC)\longrightarrow L^2(\Omega_0,\CC)$. We show that the conjugation by $h^\sharp$ provides a new 
equation of the form 
\begin{equation}\label{eq_abstract2}
i\partial_t v= h^\sharp(t)H(t)h_\sharp(t) v~~~~~\text{ on }\Omega_0\tag{$\ast\ast$}
\end{equation} 
where $h_\sharp=(h^\sharp)^{-1}$. The Hamiltonian $H(t)$ is a magnetic Laplacian 
operator of the form$$H(t)=-(\div_x+iA)\circ(\grad_x+iA)-|A|^2$$
where $A$ is an explicit magnetic potential depending on the deformation 
of the domain $\Omega(t)$. The formulation 
\eqref{eq_abstract2} enables to ensure the existence of weak and strong solutions of the initial problem 
\eqref{eq_abstract} on $\Omega(t)$ endowed with Dirichlet boundary conditions. In addition, it also indicates that 
the correct Neumann type boundary conditions for \eqref{eq_abstract} are not the homogeneous but the magnetic 
ones
$$\partial_\nu u(t)+i\la\nu| A\ra u(t)=0,$$
even though \eqref{eq_abstract} has no magnetic term. All the previous 
results are also studied in presence of diffusion coefficients as well as magnetic and electric potentials. 
Finally, we prove some associated byproducts as an adiabatic result for slow deformations of the domain and a 
time-dependent version of the so-called ``Moser's trick''. We use this outcome in order to simplify 
Equation \eqref{eq_abstract2} and  to guarantee the well-posedness for slightly less regular deformations of 
$\Omega(t)$. 

\vspace{3mm}

\end{abstract}

%%%%%%%%%%%%%%%%%%%%%%%%%%%%%%%%%%%%%%%%%%%%%%%%%%%%%%%%%%%%%%%%%%%%%%%%%%%%%%%%
%%%%%%%%%%%%%%%%%%%%%%%%%%%%%%%%%%%%%%%%%%%%%%%%%%%%%%%%%%%%%%%%%%%%%%%%%%%%%%%%

\section{Main results}\label{section_intro}
In this article, we study the well-posedness of the Schr\"odinger equation
\begin{equation}\label{SE_1}
i\ddd_t u(t, x)=-\Delta u(t, x),\ \ \ \ \ \ \ t\in I~,~~x\in \Omega(t)
\end{equation}
where $I$ is an interval of times and $t\in I\mapsto \Omega(t)\subset\RR^N$ is a time-dependent family of bounded 
domains of $\RR^N$ with $d\geq 1$. We consider the cases of Dirichlet boundary conditions and of suitable magnetic 
Neumann boundary conditions.
This kind of problem is very natural when we consider a quantum particle confined in a structure which
deforms in time. 

The Schr\"odinger equation in moving domains has been widely studied in literature and an example is the classical 
article of Doescher and Rice \cite{Doescher-Rice}. For other references on the subject, we mention \cite{BMT,Beauchard,BBL,BBFSV,DiMartino-Facchi,Knobloch-Krechetnikov,Makowski-Peplowski,Moyano,Pinder,Rouchon}. In most of these references, \eqref{SE_1} is studied in dimension $d=1$ or in higher dimensions with symmetries as the radial case or the translating case. From this perspective, the purpose of this work is natural: we aim to study the well-posedness of \eqref{SE_1} in a very general framework.

\begin{figure}[ht]
\begin{center}
\resizebox{10cm}{!}{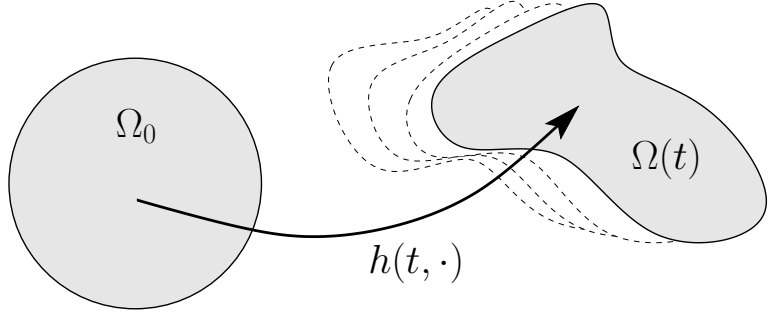}
\end{center}
\caption{\it The family of diffeomorphisms $(h(t,\cdot))_{t\in I}$ which enables to go back to a fixed domain $\Omega_0$. \label{fig-intro}} 
\end{figure}

The difficulty of considering an equation in a moving domain as \eqref{SE_1} is that the phase space 
$L^2(\Omega(t),\CC)$ and thus the operator $\Delta=\Delta(t)$ depend on the time. The usual method adopted in these kinds of problems consists of transforming $\Omega(t)$ in a bounded reference domain $\Omega_0\subset\RR^N$. Such transformation is then used in order to bring back the Schr\"odinger equation \eqref{SE_1} in an equivalent equation in the phase space 
$L^2(\Omega_0, \CC)$, which does not depend on the time. To this purpose, one can introduce a 
family of diffeomorphisms $(h(t,\cdot))_{t\in I}$ such that 
for each $t\in I$, $h(t,\cdot)$ is a $\Cc^p-$diffeomorphism from $\overline \Omega_0$ onto $\overline\Omega(t)$ 
with $p\geq 1$ (see Figure \ref{fig-intro}). Assume in addition that the function $t\in I\mapsto h(t,y)\in\RR^N$ 
is of class $\Cc^q$ with respect to the time with $q\geq 1$.

In order to bring back the Schr\"odinger equation \eqref{SE_1} in an equivalent equation in $L^2(\Omega_0, 
\CC)$, one can introduce the pullback operator  
\begin{equation}\label{pullback}
h^*(t)~:~\phi\in L^2(\Omega(t),\CC)~\longmapsto~\phi\circ h=\phi(h(t,\cdot)) \in L^2(\Omega_0,\CC)
\end{equation}
and its inverse, the pushforward operator, defined by 
\begin{equation}\label{pushforward}
h_*(t)~:~\psi\in L^2(\Omega_0,\CC)~\longmapsto~\psi\circ h^{-1}=\psi(h^{-1}(t,\cdot)) \in L^2(\Omega(t),\CC)~.
\end{equation}
If we compute the equation satisfied by $w=h^*u$ when $u$ is solution of \eqref{SE_1} at least in a formal 
sense, then we find that \eqref{SE_1} becomes
\begin{equation}\label{SE_2}
i\ddd_t w(t,y)~ = ~-
\frac{1}{|J|}\div_{y}\Big(|J| J^{-1} (J^{-1})^t \grad_y w(t,y)\Big) + i \la J^{-1}\partial_t h(t,y) | \grad_y 
w(t,y)\ra,\ \ \ \ \ \ \ \  t\in I~,~~y\in \Omega_0,
\end{equation}
where $J=J(t,y)=D_y h(t,y)$ is the Jacobian matrix of $h$, $|J|$ stands for $|\det(J)|$ and $\la \cdot|\cdot\ra$ 
corresponds to the scalar product in $\CC^N$ (see Remark \ref{transformed_dynamics_not_unitary} in Section 
\ref{section_moving} 
for the computations).

A possible way to give a sense to the Schr\"odinger equation \eqref{SE_1} consists in proving that \eqref{SE_2}, 
endowed with some boundary conditions, generates a well-posed flow in $L^2(\Omega_0,\CC)$. This can be locally done by exploiting some specific properties of the Schr\"odinger equation with perturbative terms of order 
one (see \cite{KPV,KPV2,KPV3,Linares-Ponce}). Nevertheless this method presents possible disadvantages for 
our purposes. Indeed, such approach does not provide the natural conservation of 
the $L^2-$norm and the Hamiltonian structure of the equation is lost in the sense that the new differential operator is no longer self-adjoint with respect to a natural time-independent hermitian structure. This fact represents an obstruction, not only to the proof of global existence of solutions, but also to the use of different techniques such as the adiabatic theory.

\vspace{3mm}

We are interested in studying the Schr\"odinger equation \eqref{SE_1} by preserving its Hamiltonian structure. 
From this perspective, it is natural to introduce the following unitary operator $h^\sharp(t)$ defined by
\begin{equation}\label{pullback_sharp}
h^\sharp(t)~:~\phi\in L^2(\Omega(t),\CC)~\longmapsto~\sqrt{|J(\cdot,t)|}\,(\phi\circ h)(t) \in L^2(\Omega_0,\CC)~.
\end{equation}
We also denote by $h_\sharp(t)$ its inverse 
\begin{equation}\label{pushforward_sharp}
h_\sharp(t)=(h^\sharp(t))^{-1}~:~\psi\mapsto \big( \psi/ \sqrt{|J(\cdot,t)|}\big)\circ h^{-1}~.
\end{equation}
Notice that the relation $\|h^\sharp(t)u\|_{L^2(\Omega_0)}=\|u\|_{L^2(\Omega(t))}$ enables to transport the conservative structure through the change of variables. A direct computation, provided in Section \ref{section_dirichlet}, shows that $u$ solves \eqref{SE_1} if and only if $v=h^\sharp u$ is solution of
\begin{align}
i\partial_t v(t,y)=&  -
\frac{1}{\sqrt{|J|}}\div_{y}\Bigg(|J| J^{-1} (J^{-1})^t \grad_y \Big(\frac{v(t,y)}{\sqrt{|J|}}\Big)\Bigg) 
\nonumber 
\\&+\frac i2 \frac{\partial_t\big(|J(t,y)|\big)}{|J|} v ~+ ~i {\sqrt{|J|}}\la J^{-1}\partial_t h(t,y) | \grad_y 
\frac{v(t,y)}{\sqrt{|J|}}\ra,\ \ \ \ \ \ \ \  t\in I~,~~y\in \Omega_0. \label{SE_3}
\end{align}
Written as it stands, this equation is not easy to handle. 
For instance, it is unclear whether the equation is of Hamiltonian type and how to compute its spectrum. The central argument of this paper is to show that Equation \eqref{SE_3} can be rewritten in the form
\begin{equation}\label{SE_Dirichlet}
i\partial_t v(t,y)=-h^\sharp\Big[\big(\div_x+iA_h\big)\circ\big(\grad_x+iA_h\big)+|A_h|^2\Big]h_\sharp v(t,y),\ \ \ \ \ \ ~~~t\in I,~~y\in\Omega_0,
\end{equation}
with $A_h(t,x)=-\frac 12 (h_*\partial_t h)(t,x)$. Now, the operator appearing in the last equation is the conjugate with respect to $h^\sharp$ and $h_\sharp$ of an explicit magnetic Laplacian. Thus, its Hamiltonian structure becomes obvious and some of its
properties, as the spectrum, may be easier to study. We refer to \cite{Bonnaillie,Fournais-Helffer,Helffer2,Helffer1,Helffer3,Pan} for different spectral results on magnetic Laplacian operators.

Using the unitary operator $h^\sharp$ rather than $h^*$ is natural and it was already done in the literature in some very specific frameworks in \cite{BMT,Beauchard,Moyano}. The same idea was also adopted to study quantum waveguides in the time-independent framework, where the magnetic field $A_h$ does not appear, see for instance \cite{Exner,Haag}. In \cite{Jurg1,Jurg2}, the transformation $h^*$ is used on manifolds with time-varying metrics. The authors assume $h$ preserving the volumes which yields $h^*=h^\sharp$. They obtain an operator involving a magnetic field similar to the one in \eqref{SE_Dirichlet}.
From this perspective, the relation between motion and magnetic field is not surprising. Physically, the momentum $p = mv$ of a moving particle of mass $m$, velocity $v$ and charge $q$ in a magnetic field $A$
must be replaced by $\widetilde p = mv + qA$. This corresponds to the magnetic field appearing in the equation \eqref{SE_Dirichlet} for the Galilean frames. An explicit example of the link between motion and magnetic field in our results can be seen in the boundary condition on a moving surface as in Figure \ref{fig-cylindre} below. This is also related to the notion of ``anti-convective derivative'' of Henry, see \cite{Henry}.
\vspace{5mm}

\needspace{5mm}
\noindent{\bf \underline{The Dirichlet boundar}y\underline{ condition}}\\[1mm]
Once Equation \eqref{SE_Dirichlet} is obtained, the Cauchy problem is easy to handle by using classical 
results of existence of unitary flows generated by time-dependent Hamiltonians. In our work, we refer 
to Theorem \ref{th_Kisynski} presented in the appendix. In the case of the simple Laplacian operator with 
Dirichlet boundary condition, we obtain our first main result which states the following.

\begin{theorem}\label{th_main_Dirichlet}
Let $I\subset \RR$ be an interval of times and let $\{\Omega(t)\}_{t\in I}\subset\RR^N$ with $N\in\NN^*$ be a 
family of domains. Assume that there exist a bounded reference domain $\Omega_0$ in $\RR^N$ and a family of 
diffeomorphisms $(h(t,\cdot))_{t\in I}\in\Cc^2(I\times \overline\Omega_0,\RR^N)$ such that 
$h(t,\overline\Omega_0)=\overline\Omega(t)$. 

Then, Equation \eqref{SE_Dirichlet} endowed with Dirichlet boundary conditions generates a unitary flow $\tilde 
U(t,s)$ on $L^2(\Omega_0)$ and we may define weak solutions of the Schr\"odinger equation 
\begin{equation}\label{SE_Dirichlet_2}
\left\{
\begin{array}{ll}i\ddd_t u(t, x)=-\Delta_x u(t, x),\ \ \ \ \ \ \ &t\in I~,~~x\in \Omega(t)\\
u_{|\partial\Omega(t)}\equiv 0&
\end{array}\right.
\end{equation}
by transporting this flow via $h_\sharp$ to a unitary flow $U(t,s):L^2(\Omega(s))\rightarrow L^2(\Omega(t))$.

Assume in addition that the diffeomorphisms $h$ are of class $\Cc^3$ with respect to the time and the space 
variable. Then, for any $u_0\in H^2(\Omega(t_0))\cap H^1_0(\Omega(t_0))$, the above flow defines a solution 
$u(t)=U(t,t_0)u_0$ in $\Cc^0(I,H^2(\Omega(t))\cap H^1_0(\Omega(t)))\cap \Cc^1(I,L^2(\Omega(t)))$ solving 
\eqref{SE_Dirichlet_2} in the $L^2-$sense.
\end{theorem}

Theorem \ref{th_main_Dirichlet} is consequence of the stronger result of Theorem \ref{th_Dirichlet} 
(Section \ref{section_dirichlet}) where we also include diffusion coefficients as well as magnetic and 
electric potentials. However, in this introduction, we consider the case 
of the free Laplacian to simplify the notations and to avoid further technicalities 

Also notice that the above result does not require the reference domain $\Omega_0$ to have any regularity. In 
particular, it may have corners such as a rectangle for example. Of course, all the domains $\Omega(t)$ have to be 
diffeomorphic to $\Omega_0$. Hence, we cannot create singular perturbations such as adding or removing corners and 
holes. Nevertheless, $\Omega(t)$ may typically be a family of rectangles or cylinders with different proportions.

\vspace{5mm}

\needspace{5mm}
\noindent{\bf \underline{The ma}g\underline{netic Neumann condition}}\\[1mm]
At first sight, one may naturally think to associate to the Schr\"odinger Equation \eqref{SE_1} with the homogenenous Neumann 
boundary conditions $\partial_\nu u(t,x)=0$ where $x\in\partial\Omega(t)$ and $\nu$ is the unit outward normal of $\ddd\Omega(t)$. However, these conditions cannot generate a unitary evolution, which is problematic for quantum dynamics. Simply consider the solution $u\equiv 1$, whose norm depends on the volume of $\Omega(t)$. 
Even when the volume of $\Omega(t)$ is constant, if $u$ solves \eqref{SE_1} with homogeneous Neumann boundary conditions, then the computation \eqref{calcul_L2} below shows that the evolution cannot be unitary,
except when $h_*\partial_th$ is tangent to $\partial\Omega(t)$ at all the boundary points, meaning that the shape of $\Omega(t)$ is in fact unchanged.

From this perspective, another interesting aspect of our result appears. The expression of Equation \eqref{SE_Dirichlet} indicates that the correct 
boundary conditions to consider are the ones associated with Neumann realization of the magnetic Laplacian operator, 
that are 
\begin{equation}\label{condi_neumann}
\partial_\nu (h_\sharp v)+i\la \nu|A_h \ra (h_\sharp v)=0~~\text{ on }\partial\Omega(t).
\end{equation}
If we denote by $\nu_0$ the unit outward normal of $\ddd\Omega_0$, then the last identity can be transposed in
\begin{equation}\label{condi_neumann2}
\Big\la (J^{-1})^t \nu_0 \Big| (J^{-1})^t\sqrt{|J|}\grad_y \big(\frac v{\sqrt{|J|}}\big) - \frac i2 (\partial_t h) 
v \Big\ra =0~~\text{ on }\partial\Omega_0~
\end{equation}
(see Remark \ref{BC_robin} for further details on the computations). Even though the conditions seems complicated 
on $\Omega_0$, they simply write as the classical magnetic 
Neumann boundary conditions for the original problem in $\Omega(t)$, see \eqref{SE_Robin_2} below. In particular, 
they exactly correspond to the ones of a planar wave bouncing off the moving surface, as it is clear in the example of Figure \ref{fig-cylindre}. 
We also notice that they perfectly match with the preservation of the energy since if $u(t)$ solves \eqref{SE_1} at least formally, then
\begin{align}
\partial_t \int_{\Omega(t)} |u(t,x)|^2 \dd x &= \int_{\partial\Omega(t)}\la \nu|h_*\partial_t h\ra |u(t,x)|^2 \dd x + 2 \Re \int_{\Omega(t)} \partial_t u(t,x)\overline{u}(t,x)\dd x \nonumber\\
&= \int_{\partial\Omega(t)}\la \nu|h_*\partial_th\ra |u(t,x)|^2 \dd x + 2 \Re \left( i \int_{\partial\Omega(t)} \partial_\nu u(t,x)\overline u(t,x)\dd x\right).\label{calcul_L2}
\end{align}
Once the correct boundary condition is inferred, we obtain the following result in the same way as the Dirichlet case.  
\begin{theorem}\label{th_main_Robin}
Let $I\subset \RR$ be an interval of times and let $\{\Omega(t)\}_{t\in I}\subset\RR^N$ with $N\in\NN^*$ be a 
family of domains. Assume that there exist a bounded reference domain $\Omega_0$ in $\RR^N$ of class $\Cc^1$ and a 
family of diffeomorphisms $(h(t,\cdot))_{t\in I}\in\Cc^2(I\times \overline\Omega_0,\RR^N)$ such that 
$h(t,\overline\Omega_0)=\overline\Omega(t)$. 

Then, Equation \eqref{SE_Dirichlet} endowed with the magnetic Neumann boundary conditions \eqref{condi_neumann} 
(or 
equivalently \eqref{condi_neumann2}) generates a unitary flow $\tilde U(t,s)$ on $L^2(\Omega_0)$ and we may define 
weak solutions of the Schr\"odinger equation 
\begin{equation}\label{SE_Robin_2}
\left\{
\begin{array}{ll}i\ddd_t u(t, x)=-\Delta_x u(t, x),\ \ \ \ \ \ \ &t\in I~,~~x\in \Omega(t)\\
\partial_\nu u(t, x) - \frac i2 \la \nu| h_*\partial_t h(t, x)\ra u(t, x)  = 0,&t\in I~,~~x\in \partial\Omega(t)
\end{array}\right.
\end{equation}
by transporting this flow via $h_\sharp$ to a unitary flow $U(t,s):L^2(\Omega(s))\rightarrow L^2(\Omega(t))$.

Assume in addition that the diffeomorphisms $h$ are of class $\Cc^3$ with respect to the time and the space 
variable. Then, for any $u_0\in H^2(\Omega(t_0))$ satisfying the magnetic Neumann boundary condition of 
\eqref{SE_Robin_2}, the above flow defines a solution $u(t)=U(t,t_0)u_0$ in $\Cc^0(I,H^2(\Omega(t)))\cap 
\Cc^1(I,L^2(\Omega(t)))$ solving \eqref{SE_Robin_2} in the $L^2-$sense and satisfying the magnetic Neumann boundary 
condition.
\end{theorem}

\begin{figure}[ht]
\begin{center}
\resizebox{12cm}{!}{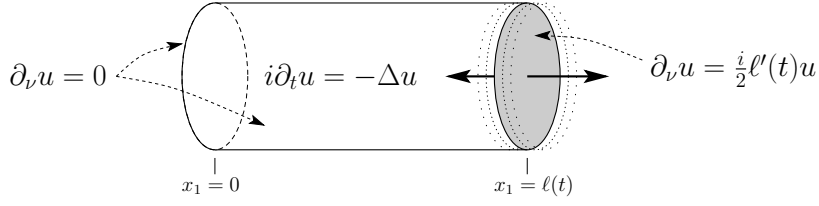}
\end{center}
\caption{\it The correct Neumann boundary conditions for the Schr\"odinger equation in a cylinder with a moving 
end. Notice the magnetic Neumann boundary condition at the moving surface, even though the equation has no 
magnetic term. See Section \ref{section_appli} for further details on the computations. \label{fig-cylindre}} 
\end{figure}

\vspace{5mm}

\needspace{5mm}
\noindent{\bf \underline{Gau}g\underline{e transformation}}\\[1mm]
As it is well known, the magnetic potential has a gauge invariance. In particular, for any $\phi$ of class $\Cc^1$ 
in space, we have
\begin{equation}\label{eq_Gauge_1}
e^{-i\phi(x)}\big[(\grad_x+iA_h)^2\big]e^{i\phi(x)}~=~ \big(\grad_x+i(A_h+\grad_x\phi)\big)^2~.
\end{equation}
Thus, it is possible to delete the magnetic term $A_h=-\frac 12 h_*\partial_t h$ by the change 
of gauge when there exists $\phi$ of class $\Cc^1$ such that 
$$\forall t\in I~,~~\forall x\in\Omega(t)~,~~(h_*\partial_t h)(t,x)=2\grad_x \phi(t,x)~.$$
In such context, the well-posedness of the equations \eqref{SE_Dirichlet_2} and \eqref{SE_Robin_2} can be 
investigated by considering 
$$w(t,y)=h^\sharp e^{-i\phi(t,x)} u:=\sqrt{|J(t,y)|}e^{-i\phi(t,h(t,y))}u(t,h(t,y)).$$
and by studying the solution of the following equation endowed with the corresponding boundary conditions
\begin{equation}\label{SE_Gauge}
i\partial_t w(t,y)=-\Big(h^\sharp\Delta_x h_\sharp +\frac 14 |\partial_t 
h(t,y)|^2-\partial_t(h^*\phi)(t,y)\Big)w(t,y),~~~t\in I~,~~y\in\Omega_0.
\end{equation}

The gauge transformation, not only simplifies Equation 
\eqref{SE_Dirichlet}, but also yields a gain of regularity in the hypotheses on $h$ adopted in the theorems 
\ref{th_main_Dirichlet} and 
\ref{th_main_Robin}. This fact follows as the main part of the new 
Hamiltonian in \eqref{SE_Gauge} does not 
contain $\partial_t h$ anymore. In details, if we consider $h\in\Cc^1_t(I,\Cc^2_x(\Omega_0,\RR^N))$, then the 
existence of weak solutions of \eqref{SE_Dirichlet_2} and \eqref{SE_Robin_2} can be guaranteed when $\phi$ is of 
class $\Cc^3$ in space and $W^{1,\infty}$ in time. The existence of strong solutions, instead, holds when $\phi$ 
is of class 
$\Cc^4$ in space and $W^{1,\infty}$ in time.

Also remark that the gauge transformation is not always possible to use. For example, if 
$\Omega(t)$ is a rotation of a square, $\partial_t h$ is not curl-free and cannot be rectified due to the presence 
of corners at which $h(t,y)$ is imposed (corners have to be send onto corners). Finally, we may also use the 
simpler gauge of the electric potential if some of the terms of \eqref{SE_Gauge} are constant, see Section 
\ref{section_appli}.

\vspace{5mm}

\needspace{5mm}
\noindent{\bf \underline{Moser's trick}}\\[1mm]
Another way to simplify Equation \eqref{SE_Dirichlet} is to use a family of diffeomorphisms $\tilde h(t)$ such 
that the determinant of the Jacobian is independent of $y$. In other words, when $\tilde J=D_y \tilde h$ satisfies 
the following identity
\begin{equation}\label{det_jac_constant}
\forall t\in I~,~~\forall y\in\Omega_0~,~~|\tilde J(t,y)|%=\frac{\meas(\Omega(t))}{\meas(\Omega_0)}
:=a(t)~. 
\end{equation}
In this case, the multiplication for the Jacobian $J$ of commutes with the spatial derivatives and then, 
Equations \eqref{SE_3} and 
\eqref{SE_Dirichlet} can be simplified in the following expression, for $A_{\tilde h}=-\frac 12 \tilde 
h_*\partial_t\tilde h$,
\begin{align}
i\partial_t v(t,y)&=-\div_y(\tilde J^{-1}(\tilde J^{-1})^t \grad_y v) + \frac i2 \frac{a'(t)}{a(t)} v + i \la 
\tilde J^{-1}\partial_t \tilde h | \grad_y v\ra \nonumber\\
&= -\tilde h^*\big[(\div_x+iA_{\tilde h})\circ(\grad_x+iA_{\tilde h})+|A_{\tilde h}|^2\big]\tilde h_* v(t,y). 
\label{SE_moser}
\end{align}
This strategy can be used, not only to simplify the equations, but also to gain some regularity since  $|\tilde 
J|$ is now constant and thus smooth. 
Therefore, $\tilde h^\sharp$ maps $H^k(\Omega(t))$ into $H^k(\Omega_0)$ as soon as $\tilde h$ is of class $\Cc^k$ 
in space.

For $t$ fixed, finding a diffeomorphism $\tilde h$ satisfying the identity \eqref{det_jac_constant} follows 
from a very famous work of 
Moser \cite{Moser}. This kind of result called ``Moser's trick'' was widely studied in literature even 
in the case of moving domains  (see Section 
\ref{section_Moser2}). Nevertheless, most of these outcomes are not interested in studying the optimal time and 
space regularity as well as their proofs are sometimes simply outlined. For this purpose, in Section 
\ref{section_Moser}, we prove the following result by following the arguments of \cite{Dacorogna-Moser}.

\begin{theorem}\label{th_moser_2}
Let $k\geq 1$, and $r\in\NN$ with $k\geq r\geq 0$. Let $\alpha\in(0,1)$ and let $\Omega_0\subset\RR^N$ be a 
connected bounded domain of class $\Cc^{k+2,\alpha}$. Let $I\subset\RR$ an interval of times and assume that there 
exists a family $(\Omega(t))_{t\in I}$ of domains such that there exists a family $(h(t))_{t\in I}$ of 
diffeomorphisms 
$$h~:~y\in \overline\Omega_0~\longrightarrow ~h(t,y)\in\overline{\Omega(t)}$$ 
which are of class $\Cc^{k,\alpha}(\overline\Omega_0,\overline{\Omega(t)})$ with respect to $y$ and of class 
$\Cc^r(\overline\Omega_0,\RR^N)$ with respect to $t$.

Then, there exists a family $(\tilde h(t))_{t\in I}$ of diffeomorphisms from $\overline\Omega_0$ onto 
$\overline{\Omega(t)}$, with the same regularity as $h$, and such that $\det(D_y\tilde h(t))$ is constant with 
respect to $y$, that is that
$$\forall y\in\Omega_0~,~~\det(D_y\tilde h(t,y))=\frac{\meas(\Omega(t))} {\meas(\Omega_0)}~.$$
\end{theorem}
Even though results similar to Theorem \ref{th_moser_2} have already been stated, the fact that $h(t)$ may simply 
be continuous with respect to the time and that $\tilde h(t)$ has the same regularity of $h(t)$ seems to be new.
There are some simple cases where we can define explicit $\tilde h$ as in dimension $N=1$ or in the examples of 
Section \ref{section_appli}, but this is not aways possible. In these last situations, Equation \eqref{SE_moser} 
may be difficult to use as $\tilde h$ is not explicit. However, Equation \eqref{SE_moser} yields a gain of 
regularity in the Cauchy problem which we resume in the following corollary.
\begin{coro}\label{coro-Moser}
Let $I\subset \RR$ be an interval of times and let $\{\Omega(t)\}_{t\in I}\subset\RR^N$ with $N\in\NN^*$ be a 
family of domains. Assume that there exist a bounded reference domain $\Omega_0$ in $\RR^N$ of class 
$\Cc^{4,\alpha}$ with $\alpha\in (0,1)$ and a family of diffeomorphisms $(h(t,\cdot))_{t\in I}\in\Cc^2_t(I,\Cc^1_x 
(\Omega_0,\RR^N))$ such that $h(t,\overline\Omega_0)=\overline\Omega(t)$. 

Then, we may define as in Theorems \ref{th_main_Dirichlet} and \ref{th_main_Robin} weak solutions of the above 
equations \eqref{SE_Dirichlet_2} and \eqref{SE_Robin_2} by considering $v=\tilde h_\sharp u$, with $\tilde h$ as 
in 
Theorem \ref{th_moser_2}, which solves the Schr\"odinger equation \eqref{SE_moser} with the corresponding boundary 
conditions and by transporting the flow of \eqref{SE_moser} via the above change of variable.

Assume in addition that $h$ belongs to $\Cc^3_t(I,\Cc^2_y(\Omega_0,\RR^N))$. Then, for any $u_0\in 
H^2(\Omega(t_0))\cap H^1_0(\Omega(t_0))$, the above flow defines a solution $u(t)$ in the space 
$\Cc^0(I,H^2(\Omega(t))\cap H^1_0(\Omega(t)))\cap \Cc^1(I,L^2(\Omega(t)))$ solving \eqref{SE_Dirichlet_2} or 
\eqref{SE_Robin_2} in the $L^2-$sense.
\end{coro}
Corollary \ref{coro-Moser} follows from the same arguments leading to the theorems 
\ref{th_main_Dirichlet} or \ref{th_main_Robin} (see Section \ref{section_cauchy}). The only difference is the gain 
of regularity in space. Indeed, the term $|J(t,y)|$ appearing in $h^\sharp$ or $h_\sharp$ is replaced 
by $|\tilde J(t,y)|=\frac{\meas(\Omega(t))}{\meas(\Omega_0)}=a(t)$ which is constant in space and then smooth. To 
this end, we simply have to replace the 
first family of diffeomorphisms $h(t)$ by the one given by Theorem \ref{th_moser_2}.

Notice that, in Corollary \ref{coro-Moser}, we have to assume that the reference domain $\Omega_0$ is smooth. If 
$\Omega(t)$ are simply of class $\Cc^1$ or $\Cc^2$, this is not a real restriction since we may choose a smooth 
reference domain and a not so smooth diffeomorphism $h$. In the case where $\Omega(t)$ has corner, as rectangles 
for example, then Corollary \ref{coro-Moser} do not formally apply. However, in the case of moving rectangles, 
finding a family of explicit diffeomorphisms $h(t)$ satisfying \eqref{det_jac_constant} is easy and the arguments 
behind Corollary \ref{coro-Moser} can be directly used, see the computations of Section \ref{section_appli}.

\vspace{5mm}

\begin{figure}[ht]
\begin{center}
\begin{center}
\resizebox{12cm}{!}{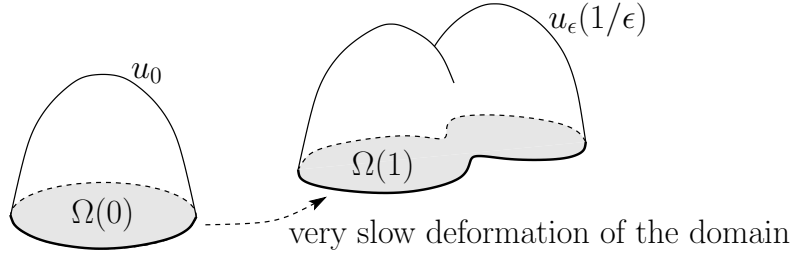}
\end{center}
\end{center}
\caption{\it The ground state of a domain $\Omega(0)$ can be transformed into the ground state of a domain 
$\Omega(1)$ if we slowly and smoothly deform $\Omega(0)$ to $\Omega(1)$ while the quantum state evolves following 
Schr\"odinger equation.  \label{fig-adiab}} 
\end{figure}

\needspace{5mm}
\noindent{\bf \underline{An exam}p\underline{le o}f\underline{ a}pp\underline{lication: an adiabatic 
result}}\\[1mm]
Consider a family of domains $\{\Omega(\tau)\}_{\tau\in [0,1]}$ of $\RR^N$ such that there exist a bounded 
reference domain $\Omega_0\subset\RR^N$ and a family of diffeomorphisms $(h(\tau,\cdot))_{\tau\in [0,1]}\in 
\Cc^3([0,1]\times\overline\Omega_0,\RR^N)$ such that $h(\tau,\overline \Omega_0)=\overline\Omega(\tau)$. Assume 
that for each $\tau\in [0,1]$, the Dirichlet Laplacian operator $\Delta$ on $\Omega(\tau)$ has a simple isolated 
eigenvalue $\lambda(\tau)$ with normalized eigenfunction $\varphi(\tau)$,  associated with a spectral projector 
$P(\tau)$, all three depending continuously on $\tau$. Following the classical adiabatic principle, we expect that 
if we start with a quantum state in $\Omega(0)$ close to $\varphi(0)$ and we deform very slowly the domain to the 
shape $\Omega(1)$, then the final quantum state should be close to $\varphi(1)$ up to a phase shift (see Figure 
\ref{fig-adiab}). The slowness of the deformation is represented by a parameter $\epsilon>0$ and we consider 
deformations between the times $t=0$ and $t=1/\epsilon$, that is the following Schr\"odinger equation
\begin{equation}\label{SE_adiabatic}
\left\{
\begin{array}{ll}i\ddd_t u_\epsilon(t, x)=-\Delta u_\epsilon(t, x),\ \ \ \ &t \in [0,1/\epsilon]~,~~x\in 
\Omega(\epsilon t)\\
u_{\epsilon}(t) \equiv 0,&\text{ on }\partial\Omega(\epsilon t) \\
u_\epsilon(t=0)= u_0 \in L^2(\Omega(0)) &
\end{array}\right.
\end{equation}
Due to Theorem \ref{th_main_Dirichlet}, we know how to define a solution of \eqref{SE_adiabatic}. Our main 
equation 
\eqref{SE_Dirichlet} provides the Hamiltonian structure required for applying the adiabatic theory, in 
contrast to the case of Equation \eqref{SE_2}. It also indicates that we should not consider the Laplacian 
operators $\Delta$ and its spectrum, but rather the magnetic Laplacian operators of Equation \eqref{SE_Dirichlet}. 
With these hints, it is not difficult to adapt the classical adiabatic methods to obtain the following result (see 
Section \ref{section_appli}).
\begin{coro}\label{coro_adiabatic}
Consider the above framework. Then, we have
$$\la P(1)u_\epsilon(1/\epsilon)|u_\epsilon(1/\epsilon)\ra ~~\xrightarrow[~~\epsilon\longrightarrow 0~~]{} ~~\la 
P(0)u_0|u_0\ra~.$$
\end{coro}
 
\vspace{5mm}

\noindent{\bf Acknowledgements:} the present work has been conceived in the vibrant atmosphere of 
the Workshop-Summer School of Benasque (Spain) ``VIII Partial differential equations, optimal design and 
numerics''. The authors would like to thank Yves Colin de Verdi\`ere and Andrea Seppi for the fruitful 
discussions on the geometric aspects of the work. They are also grateful to G\'erard Besson and Emmanuel Russ for 
the suggestions on the proof of the Moser's trick and to Alain Joye for the advice on the adiabatic theorem. This 
work was financially supported by the project {\it 
ISDEEC} ANR-16-CE40-0013.

%%%%%%%%%%%%%%%%%%%%%%%%%%%%%%%%%%%%%%%%%%%%%%%%%%%%%%%%%%%%%%%%%%%%%%%%%%%%%%%%
%%%%%%%%%%%%%%%%%%%%%%%%%%%%%%%%%%%%%%%%%%%%%%%%%%%%%%%%%%%%%%%%%%%%%%%%%%%%%%%%

\section{Moving domains and change of variables}\label{section_moving}
The study of the influence of the domain shape on a PDE problem has a long history, in particular in shape 
optimization. The interested reader may consider \cite{Henry} or \cite{Murat-Simon} for example. The 
basic strategy of these kinds of works is to bring back the problem in a fixed reference domain $\Omega_0$ via 
diffeomorphisms. We thus have to compute the new differential operators in the new coordinates. The formulae 
presented in this section are well known (see \cite{Henry} for example). We recall them for sake of completeness 
and also to fix the notations. 

\vspace{3mm}

Let $\Omega_0$ be a reference domain and let $h(t)$ be a $\Cc^2-$diffeorphism mapping $\Omega_0$ onto $\Omega(t)$ 
(see Figure \ref{fig-intro}). In this section, we only work with fixed $t$. For this reason, we omit the time 
dependence when it is not necessary and forgot the question of regularity with respect to the time.

From now on, we 
denote by $x$ the points in $\Omega(t)$ and by $y$ the ones in $\Omega_0$. For any matrix $A$, $|A|$ stands for 
$|\det(A)|$ and $\la\cdot|\cdot\ra$ is the scalar product 
in $\CC^N$ with the convention 
$$\la v|w\ra=\sum_{k=1}^N\overline{v_k} \,w_k~, \ \ \ \ \ \  \ \ \ \forall v,w\in \CC^N.$$
We use the pullback and pushforward 
operators defined by \eqref{pullback} and \eqref{pushforward}.
We denote by $J=J(t,y)=D_y h(t,y)$ the Jacobian matrix of $h$. The basic rules of differential calculus give 
that
\begin{equation}\label{jacobian_inverse}
h^*D_x(h^{-1})=(D_yh)^{-1}:=J^{-1}~,\ \ \ \ \ \ \ \ |J^{-1}|=|J|^{-1}
\end{equation}
and 
\begin{equation}\label{change_variable_x}
\int_{\Omega(t)} f(x) \dd x~=~\int_{\Omega_0} |J(t,y)| (h^*f)(y)\dd y~.
\end{equation}
The result below is due to the chain rule which implies the following identity
\begin{equation}\label{deriv_partiel_y}
\ddd_{y_j}(h^*f)=\la \ddd_{y_j}h|h^*\grad_{x}f\ra.
\end{equation}

\begin{prop}\label{prop_grad}
For every $g\in H^1(\Omega_0,\CC)$, 
$$(h^*\grad_x h_*)g(y)=\big(J(t,y)^{-1}\big)^t\cdot\grad_y g(y)$$
where $\big(J^{-1}\big)^{t}$ is the transposed matrix of $(D_yh)^{-1}$.
\end{prop}
\begin{demo}
We apply \eqref{deriv_partiel_y} with $f=h_*g$ in order to obtain the identity $\grad_y g=J(t,y)^t h^*(\grad_x f)$.
\end{demo}

By duality, we obtain the corresponding property for the divergence.
\begin{prop}\label{prop_div}
For every vector field $A\in H^1(\Omega_0,\CC^N)$,
$$(h^*\div_{x}h_*)A(y)=\frac{1}{|J(t,y)|}\div_{y}\big( |J(t,y)|J(y,t)^{-1} A(y)\big).$$
Let $\nu$ and $\sigma$ respectively be the unit outward normal and the measure on the boundary of $\Omega(t)$. 
Let $\nu_0$ and $\sigma_0$ respectively be the unit outward normal and the measure on the boundary of $\Omega_0$. 
For every $A\in H^1(\Omega_0,\CC^N)$ and $g\in H^1(\Omega_0,\CC)$, 
$$\int_{\ddd \Omega(t)} \la \nu |(h_*A)(x) \ra (h_*g)(x)\dd\sigma=\int_{\ddd \Omega_0}|J(t,y)|\big\la 
\nu_0\big | J(t,y)^{-1}A(y)\big\ra g(y) \dd\sigma_0.$$
\end{prop}
\begin{demo} Let $\varphi\in\Cc^\infty_0(\Omega_0,\CC)$ be a test function and let $A\in H^1(\Omega_0,\CC^N)$. 
Applying the divergence theorem and the above formulas, we obtain
\begin{align*}
\int_{\Omega_0} |J(t,y)|(h^*\div_x h_*A)(y)\varphi(y)\dd y &=
\int_{\Omega(t)} (\div_x (h_*A))(x)(h_*\varphi)(x)\dd x\\
 &= -\int_{\Omega(t)} \la h_*\overline A(x)| \grad_{x}(h_*\varphi)(x)\ra\dd x\\
&=-\int_{\Omega_0} |J(y,t)|\big\la \overline A(y)| h^*(\grad_{x}h_*\varphi)(y)\big\ra\dd y\\
&=-\int_{\Omega_0} |J(y,t)|\big\la \overline A(y)| (J(t,y)^{-1})^t \grad_y \varphi(y)\big\ra\dd y\\
&=-\int_{\Omega_0} \big\la |J(y,t)|J(t,y)^{-1}\overline A(y)| \grad_y \varphi(y)\big\ra\dd y\\
&=\int_{\Omega_0} \div_y\big[|J(y,t)|J(t,y)^{-1}A(y)\big]\varphi(y)\dd y.
\end{align*}
The first statement follows from the density of $\Cc^\infty_0(\Omega_0,\CC)$ in $L^2(\Omega_0)$. 
The second one instead is proved by considering the border terms appearing when we proceed as above with $A\in 
H^1(\Omega(t),\CC^N)$ and $g\in 
H^1(\Omega(t),\CC)$. In this context, we use twice the divergence theorem and we obtain that
\begin{align*}
\int_{\Omega_0} |J(t,y)|(h^*\div_x h_*A)(y)g(y)\dd y = &
\int_{\partial\Omega(t)} \la \nu|(h_*A)(x) \ra (h_*g)(x)\dd\sigma \\
&-\int_{\ddd \Omega_0}|J(t,y)|\big\la \nu_0\big| J(t,y)^{-1}A(y)\big\ra g(y) \dd\sigma_0\\
&+ \int_{\Omega_0} \div_y\big(|J(y,t)|J(t,y)^{-1}A(y)\big)g(y)\dd y.
\end{align*}
The last relation yields the equality of the boundary terms since the first statement is valid in the $L^2$ 
sense.
\end{demo}

By a direct application of the previous propositions, we obtain the operator associated with 
$\Delta_x=\div_x(\grad_x\cdot)$.
\begin{prop}\label{prop_Lapl}
For every $g\in H^2(\Omega_0,\CC)$,  
$$(h^*\Delta_xh_*)g(y)=\frac{1}{|J(t,y)|}\div_{y}\Big(|J(t,y)| J(t,y)^{-1} (J(t,y)^{-1})^t \grad_y g(y)\Big).$$ 
\end{prop}

\begin{remark}\label{transformed_dynamics_not_unitary}
For an illustration, let us check, at least formally, that $u(t,x)$ solves the Schr\"odinger equation \eqref{SE_1} 
in $\Omega(t)$ if and only if $w=h^*u$ solves \eqref{SE_2}. We have
\begin{align*}
i\partial_t w(t,y)&= i\partial_t \big(u(t,h(t,y))\big)\\
&= (i\partial_t u)(t,h(t,y)) + i \la \partial_t h(t,y) |(\grad_x u)(t,h(t,y))\ra\\
&=(-\Delta_x u)(t,h(t,y)) + i \la \partial_t h(t,y) |(\grad_x u)(t,h(t,y))\ra\\
&=-(h^*\Delta_x h_*w)(t,y) + i \la\partial_t h(t,y)|(h^*\grad_x h_*w)(t,y)\ra\\
&=-\frac{1}{|J|}\div_{y}\Big(|J| J^{-1} (J^{-1})^t \grad_y g(t,y)\Big) + i \la \partial_t h(t,y) 
|(J^{-1})^t\grad_y 
w(t,y)\ra\\
&=-\frac{1}{|J|}\div_{y}\Big(|J| J^{-1} (J^{-1})^t \grad_y g(t,y)\Big) + i \la J^{-1}\partial_t h(t,y) | \grad_y 
w(t,y)\ra.
\end{align*}

\end{remark}
%%%%%%%%%%%%%%%%%%%%%%%%%%%%%%%%%%%%%%%%%%%%%%%%%%%%%%%%%%%%%%%%%%%%%%%%%%%%%%%%
%%%%%%%%%%%%%%%%%%%%%%%%%%%%%%%%%%%%%%%%%%%%%%%%%%%%%%%%%%%%%%%%%%%%%%%%%%%%%%%%
\section{Main results: The Cauchy problem}\label{section_cauchy}
\subsection{Proof of Theorem \ref{th_main_Dirichlet}: the Dirichlet case}\label{section_dirichlet}
Let $I$ be a time interval and let $\Omega(t)=h(t,\Omega_0)$ be a family of moving domains.
In this section, we consider a second order Hamiltonian operator of the type
\begin{align}
H(t)&=-\Big[D(t,x)\grad_x+iA(t,x)\Big]^2+V(t,x) \label{eqH}
\end{align}
where
\begin{align}
\Big[D(t,x)\grad_x+iA(t,x)\Big]^2:=\Big(\div_x \big(D(t,x)^t\, \cdot\,) ,+\,i\la  
A(t,x)|\ \cdot \ \ra\Big)\circ\Big(D(t,x)\grad_x\,+\,iA(t,x)\Big)\label{eqH_notation}
\end{align}
and
\begin{itemize}
\item $D\in\Cc^2_t(I,\Cc^1_x(\RR^N,\Mc_N(\RR)))$ are symmetric diffusions coefficients such that there exists 
$\alpha>0$ satisfying
$$\forall t\in I~,~~\forall x\in\RR^N~,~~\forall \xi\in\CC^N~,~~\la D(x,t)\xi | D(x,t)\xi \ra \geq \alpha 
|\xi|^2~;$$
\item $A\in\Cc^2_t(I,\Cc^1_x(\RR^N,\RR^N))$ is a magnetic potential;
\item $V\in\Cc^2_t(I,L^\infty_x(\RR^N,\RR))$ is an electric potential.% such that there exists $\beta\in\RR$ such 
%that 
%$$\forall t\in I~,~~\forall x\in\RR^N~,~~V(t,x)\geq \beta;$$
\end{itemize}
In this subsection, we assume that $H(t)$ is associated with Dirichlet boundary conditions and then its domain is 
$H^2(\Omega(t),\CC)\cap H^1_0(\Omega(t),\CC)$. Of course, the typical example of such Hamiltonian is the Dirichlet 
Laplacian $H(t)=\Delta_x$ defined on 
$H^2(\Omega(t),\CC)\cap H^1_0(\Omega(t),\CC)$. We consider the equation
\begin{equation}\label{eq_dirichlet_1}
\left\{\begin{array}{l}
i\partial_t u(t,x) =H(t)u(t,x),~~~~t\in I~,~x\in \Omega(t)\\
u_{|\partial\Omega(t)}\equiv 0 
\end{array}\right. 
\end{equation}
We use the pullback and pushforward operators $h^\sharp$ and $h_\sharp$ defined by \eqref{pullback_sharp} and 
\eqref{pushforward_sharp}. We also use the notations of Section \ref{section_moving}. Notice that $h^\sharp$ is an 
isometry from $L^2(\Omega(t))$ onto $L^2(\Omega_0)$. Moreover, if $h$ is a class $\Cc^2$ with respect to the 
space, 
then $h^\sharp$ is an isomorphism from $H^1_0(\Omega(t))$ onto $H^1_0(\Omega_0)$. We set
\begin{equation}\label{def_v}
v(t,y)~:=~(h^\sharp u)(t,y)~=~\sqrt{|J(t,y)|}\,u(t,h(t,y))~. 
\end{equation}
At least formally, if $u$ is solution of \eqref{eq_dirichlet_1}, then a direct computation shows that 
\begin{align}
i\partial_t v(t,y) & =  h^\sharp (i\partial_t u)(t,y) + \frac i2 
\frac{\partial_t\big(|J(t,y)|\big)}{\sqrt{|J(t,y)|}} u(t,y) + i \sqrt{|J(t,y)|} \big\la \partial_t h(t,y) |  
h^*(\grad_x u) \big\ra \nonumber\\
& =  \big(h^\sharp H(t) h_\sharp\big)v(t,y) + \frac i2 \frac{\partial_t\big(|J(t,y)|\big)}{|J(t,y)|} v(t,y) + i 
\big\la \partial_t h(t,y) | (h^\sharp \grad_x h_\sharp) v(t,y) \big\ra  \nonumber\\
& =  \big(h^\sharp H(t) h_\sharp\big)v(t,y) + \frac i2 \frac{\partial_t\big(|J(t,y)|\big)}{|J(t,y)|} v(t,y)  
+ h^\sharp\Big[ i \big\la (h_*\partial_t h)(t,x) | \grad_x \cdot \big\ra \Big] h_\sharp 
v(t,y).\label{eq_dirichlet_2}
\end{align}
The first operator $(h^\sharp H(t) h_\sharp)$ is a simple transport of the original operator $H(t)$ and it is 
clearly self-adjoint. Both other terms from the last relation come instead from the time derivative of 
$h^\sharp$. Since $h^\sharp$ is unitary, we expect that, their sum is also a self-adjoint operator. 
Notice that the last term may be expressed explicitly by Proposition \ref{prop_grad} to obtain Equation 
\eqref{SE_3}. However, we would like to keep the conjugated form to obtain Equation \eqref{SE_Dirichlet}. Due to 
Proposition \ref{prop_jacobi} in the appendix and 
\eqref{jacobian_inverse}, we have
\begin{align*}
\frac{\partial_t\big(|J(t,y)|\big)}{|J(t,y)|} & = \Tr\Big(J(t,y)^{-1}\cdot\partial_t  J(t,y)\Big)= 
\Tr\Big(\partial_t J(t,y)\cdot J(t,y)^{-1}\Big) \\
&= h^* \Tr\Big(\big(\partial_t D_y h\big) \circ h^{-1}\cdot(D_y h)^{-1}\circ h^{-1}\Big)\\
& = h^* \Tr\Big(\big( D_y \partial_t h\big) \circ h^{-1}\cdot D_x (h^{-1})\Big)=  h^* \Tr\Big(D_x \big( 
(\partial_t h)\circ h^{-1}\big) \Big)\\
&= h^* \div_x \big(h_*(\partial_t h)\big)~.
\end{align*}
Since $h^\sharp$ differs from $h^*$ by a multiplication, the conjugacy of a multiplicative operation by either 
$h^*$ or $h^\sharp$ gives the same result. We obtain that 
\begin{align}
\frac i2 \frac{\partial_t\big(|J(t,y)|\big)}{|J(t,y)|} v(t,y) &= \frac i2  \Big[h^* \div_x \big(h_*(\partial_t 
h)\big)\Big] v = \frac i2 h^*\Big[\div_x \big(h_*(\partial_t h)\big) h_*v \Big]\nonumber\\
&=  \frac i2 h^\sharp\Big[\div_x \big(h_*(\partial_t h)\big) h_\sharp v \Big].\label{derivative_jacobian}  
\end{align}
We combine the result of \eqref{derivative_jacobian} and half of the last term of \eqref{eq_dirichlet_2} by using 
the chain rule
$$\frac i2 \div_x \big(h_*(\partial_t h)\big) u +  i \big\la (h_*\partial_t h)| \grad_x u \big\ra =  \frac i2 
\div_x \Big( \big(h_*(\partial_t h)\big)u\Big) +  \frac i2 \big\la (h_*\partial_t h)| \grad_x u \big\ra.
$$
We finally obtain
\begin{equation}\label{eq_dirichlet_3}
i\partial_t v(t,y)~=~h^\sharp\Big[H(t) \cdot  ~-~ i \div_x \big(A_h(t,x) \cdot \big) ~-~ i \big\la A_h(t,x)\big| 
\grad_x \cdot \big\ra\Big] h_\sharp v(t,y)
\end{equation}
where $A_h(t,x)=-\frac 12 (h_*\partial_t h)(t,x)$ is a magnetic field generated by the change of referential.
The whole operator of \eqref{eq_dirichlet_3} can be seen as a modification of the magnetic and electric terms of 
$H(t)$. Indeed, using \eqref{eqH}, the equation \eqref{eq_dirichlet_3} becomes
\begin{equation}\label{eq_dirichlet_4}
\left\{
\begin{array}{ll}i\partial_t v(t,y)~=~\big(h^\sharp\tilde H_h(t)h_\sharp\big) v(t,y),~~~&t\in I~,~~y\in\Omega_0\\
v_{|\partial \Omega_0}\equiv 0
\end{array}\right.
\end{equation}
where, using the same standard notation of \eqref{eqH_notation},
\begin{equation}
\tilde H_h(t)=-\Big[D(t,x)\grad_x+i\tilde A_h(t,x)\Big]^2+\tilde V_h(t,x) \label{eqH2}
\end{equation}
with 
$$\tilde A_h(t,x)=A(t,x)+\big(D(t,x)^{-1}\big)^t\cdot A_h(t,x)~~,~~~~~~~~~A_h=-\frac 12 h_*(\partial_t h),$$ 
$$\tilde V_h(t,x)=V(t,x)-\big|\big( D(t,x)^{-1}\big)^t\cdot A_h(t,x)\big|^2-2\big\la \big(D(t,x)^{-1}\big)^t\cdot 
A_h(t,x)\big|A(t,x) \big\ra~.$$
Notice that, in the simplest case of $H(t)$ being the Dirichlet Laplacian operator, we obtain Equation 
\eqref{SE_Dirichlet} discussed in the introduction.  Similar computations are provided in \cite{Jurg1,Jurg2}.

\vspace{3mm}

It is now possible to prove the main result of this section which generalizes Theorem \ref{th_main_Dirichlet} of 
the introduction.

\begin{theorem}\label{th_Dirichlet}
Assume that $\Omega_0\subset\RR^N$ is a bounded domain (possibly irregular). Let $I\subset \RR$ be an interval of 
times. Assume that $h(t):\overline\Omega_0\longrightarrow \overline\Omega(t):=h(t,\overline\Omega_0)$ is a family 
of diffeomorphisms which is of class $\Cc^2$ with respect to the time and the space variable. Let $H(t)$ be 
an operator of the form \eqref{eqH} with domain $H^2(\Omega(t),\CC)\cap H^1_0(\Omega(t),\CC)$.

Then, Equation \eqref{eq_dirichlet_4} generates a unitary flow $\tilde U(t,s)$ on $L^2(\Omega_0)$ and we may 
define 
weak solutions of the Schr\"odinger equation \eqref{eq_dirichlet_1} on the moving domain $\Omega(t)$ by 
transporting this flow via $h_\sharp$, that is setting $U(t,s)=h_\sharp \tilde U(t,s) h^\sharp$.

Assume in addition that the diffeomorphisms $h$ are of class $\Cc^3$ with respect to the time and the space 
variable. Then, for any $u_0\in H^2(\Omega(t_0))\cap H^1_0(\Omega(t_0))$, the above flow defines a solution 
$u(t)=U(t,t_0)u_0$ in $\Cc^0(I,H^2(\Omega(t))\cap H^1_0(\Omega(t)))\cap \Cc^1(I,L^2(\Omega(t)))$ solving 
\eqref{eq_dirichlet_1} in the $L^2-$sense.
\end{theorem}
\begin{demo}
Assume that $I$ is compact, otherwise it is sufficient to cover $I$ with compact intervals and to glue the 
unitary flows defined on each one of them. We simply apply Theorem \ref{th_Kisynski} in 
appendix. We notice $A_h=-\frac 12 h_*(\partial_th)$ is a bounded term in 
$\Cc^1_t(I,\Cc^2_x(\overline\Omega(t),\RR))$. 
First, the assumptions on $H(t)$ and the regularity of the diffeomorphisms $h(t,\cdot)$ imply that 
$\tilde H_h(t)$ defined by \eqref{eqH2} is a well defined self-adjoint operator on $L^2(\Omega(t))$ with domain 
$H^2(\Omega(t))\cap H^1_0(\Omega(t))$. Second, it is of class $\Cc^1$ with respect to the time and, for every 
$u\in H^2\cap H^1(\Omega(t))$,
\begin{align*}
\big\la\tilde H_h(t)u|u\big\ra_{L^2(\Omega(t))}&= \int_{\Omega(t)} \big|\big(D(t,x)\grad_x +i \tilde A_h\big)u 
(x)\big|^2 + \tilde V_h(t,x)|u(x)|^2\dd x\\
&\geq \gamma \|u\|_{H^1(\Omega(t))}^2 - \kappa \|u\|_{L^2(\Omega(t))}^2
\end{align*}
for some $\gamma>0$ and $\kappa\in\RR$. 
Let $\Hc(t)=h^\sharp(t)\circ \tilde H_h(t)\circ h_\sharp(t)$. Since $h^\sharp$ is an isometry from 
$L^2(\Omega(t))$ 
onto $L^2(\Omega_0)$ continuously mapping $H^1(\Omega(t))$ into $H^1(\Omega_0)$, the above properties of $\tilde 
H_h(t)$ are also valid for $\Hc(t)$ because for every $v\in H^2(\Omega_0)\cap H^1_0(\Omega_0)$
$$ \big\la \Hc(t)v|v\big\ra_{L^2(\Omega_0)} ~=~\big\la h^\sharp \tilde H_h(t) h_\sharp 
v|v\big\ra_{L^2(\Omega_0)}~=~ \big\la\tilde H_h(t)h_\sharp v|h_\sharp v\big\ra_{L^2(\Omega(t))}~.$$
The only problem concerns the regularities. Indeed, the domain of $\Hc(t)$ is 
$$D(\Hc(t))~=~\{v\in L^2(\Omega_0)~,~h_\sharp(t)v\in H^2(\Omega(t))\cap H^1_0(\Omega(t))\}$$
and if $h$ is only $\Cc^2$ with respect to the space, then $D(\Hc(t))$ is not necessarily $H^2(\Omega_0)\cap 
H^1_0(\Omega_0)$ due to the presence of the Jacobian of $h$ in the definition \eqref{pushforward_sharp} of 
$h_\sharp$. However, if $h$ is of class $\Cc^2$, then $h_\sharp$ transports the $\Cc^1-$regularity in space as well 
as the boundary condition and 
$$D(\Hc(t)^{1/2})~=~\{v\in L^2(\Omega_0)~,~h_\sharp(t)v\in H^1_0(\Omega(t))\}~=~H^1_0(\Omega_0)$$
does not depend on the time. We can apply Theorem \ref{th_Kisynski} in the appendix and obtain the unitary flow 
$\tilde U(t,s)$. The flow $U(t,s)=h_\sharp \tilde U(t,s) h^\sharp$ on $L^2(\Omega(t))$ is then well defined but 
corresponds to solutions of \eqref{eq_dirichlet_1} only in a formal way.

Assume finally that the diffeomorphisms $h$ are of class $\Cc^3$ with respect to the time and the space variable. 
Then, we have no more problems with the domains and 
$$D(\Hc(t))~=~\{v\in L^2(\Omega_0)~,~h_\sharp(t)v\in H^2(\Omega(t))\cap H^1_0(\Omega(t))\}~=~H^2(\Omega_0)\cap 
H^1_0(\Omega_0)$$
which does not depend on the time. Since $H(t)$ is now of class $\Cc^2$ with respect to the time, we may apply the 
second part of Theorem \ref{th_Kisynski} and obtain strong solutions of the equation on $\Omega_0$, which are 
transported to strong solutions of the equation on $\Omega(t)$.
\end{demo}

%%%%%%%%%%%%%%%%%%%%%%%%%%%%%%%%%%%%%%%%%%%%%%%%%%%%%%%%%%%%%%%%%%%%%%%%%%%%%%%%
\subsection{Proof of Theorem \ref{th_main_Robin}: the magnetic Neumann case}
In the previous section, the well-posedness of the 
dynamics of \eqref{eq_dirichlet_4} is ensured by studying the self-adjoint operator 
$\tilde H_h(t)$ defined in \eqref{eqH2} and with domain $H^2(\Omega(t))\cap H^1_0(\Omega(t))$. This operator 
corresponds to the following quadratic form in $H^1_0(\Omega(t))$
$$q(\psi):=\int_{\Omega(t)}\Big|D(t,x)\grad\psi(x)+i\tilde 
A_h(t,x)\psi(x)\Big|^2\dd x+\int_{\Omega(t)}\tilde V_h(t,x)|\psi(x)|^2\dd x,\ \ \ \ \ \forall \psi\in 
H^1_0(\Omega(t))$$
which is the Friedrichs extension of the quadratic form $q$ defined on $\Cc^\infty_0(\Omega(t))$. Now, it is 
natural to consider the Friedrichs extension of $q$ defined in $H^1(\Omega(t))$.
This corresponds to the Neumann realization of the magnetic Laplacian $\tilde H_h(t)$ defined by \eqref{eqH2} and 
with domain
\begin{equation}\label{domain_robin} 
D\big(\tilde H_h(t)\big)= \Big\{u\in H^2(\Omega(t),\CC)\ :\ \big\la\,\nu\,\big|\,D\grad_x u + i\tilde A_h u\big\ra 
(t,x)  = 0,\ \ \forall x\in \ddd\Omega(t)\Big\}
\end{equation}
(see \cite{Fournais-Helffer} for example).
Such as the well-posedness of the Schr\"odinger equation on moving domains can be achieve when it is endowed with 
Dirichlet boundary conditions, the same result can be addressed in this new framework. Indeed, the operators 
$\tilde H_h(t)$ endowed with the domain \eqref{domain_robin} are still self-adjoint and bounded from below. 
The arguments developed in the previous section lead to the well-posedness in $\Omega_0$ of the 
equation
\begin{equation}\label{eq_robin_2}
\left\{
\begin{array}{ll}i\partial_t v(t,y)~=~\big(h^\sharp\tilde H_h(t)h_\sharp\big) v(t,y),~~~&t\in I~,~~y\in\Omega_0\\
\Big(h^\sharp \big\la\,\nu\,\big|D\grad_x (h_\sharp v)+ i\tilde A_h (h_\sharp v) \big\ra\Big) (t,y)  = 
0,~~~~~&t\in 
I~,~x\in \ddd\Omega_0
\end{array}\right.
\end{equation}
Going back to the original moving domain $\Omega(t)$, \eqref{eq_robin_2} becomes
\begin{equation}\label{eq_robin_1}
\begin{split}\begin{cases}
i\partial_t u(t,x) =H(t)u(t,x),~~~~~~~~~~&t\in I~,~x\in \Omega(t)\\
\big\la\,\nu\,\big|D\grad_x u+  i\tilde A_h u\big\ra 
(t,x)  = 0,~~~~~&t\in I~,~x\in \ddd\Omega(t).
\end{cases}\end{split}
\end{equation}
It is noteworthy that in this case, the modified magnetic potential $\tilde A_h$ appears in the original equation. 
In particular, notice that the boundary condition of \eqref{eq_robin_1} is not the natural one associated with $H(t)$ 
for $t$ fixed. 

We resume in the following theorem the well-posedness  of the dynamics when the Neumann magnetic boundary 
conditions are considered. The result is a 
generalization of Theorem \ref{th_main_Robin} in the introduction. 
\begin{theorem}\label{th_robin}
Assume that $\Omega_0\subset\RR^N$ is a bounded domain (possibly irregular). Let $I\subset \RR$ be an interval of 
times. Assume that $h(t):\overline\Omega_0\longrightarrow \overline\Omega(t):=h(t,\overline\Omega_0)$ is a family 
of diffeomorphisms which is of class $\Cc^2$ with respect to the time and the space variable. Let $H(t)$ be a 
of the form \eqref{eqH} and with domain \eqref{domain_robin}.

Then, Equation \eqref{eq_robin_2} generates a unitary flow $\tilde U(t,s)$ on $L^2(\Omega_0)$ and we may define 
weak solutions of the Schr\"odinger equation \eqref{eq_robin_1} on the moving domain $\Omega(t)$ by 
transporting this flow via $h_\sharp$, that is setting $U(t,s)=h_\sharp \tilde U(t,s) h^\sharp$.

Assume in addition that the diffeomorphisms $h$ are of class $\Cc^3$ with respect to the time and the space 
variable. Then, for any $u_0\in H^2(\Omega(t_0))\cap H^1_0(\Omega(t_0))$, the above flow defines a solution 
$u(t)=U(t,t_0)u_0$ in $\Cc^0(I,H^2(\Omega(t))\cap H^1_0(\Omega(t)))\cap \Cc^1(I,L^2(\Omega(t)))$ solving 
\eqref{eq_robin_1} in the $L^2-$sense.
\end{theorem}
\begin{demo} {\it Mutatis mutandis}, the proof is the same as the one of Theorem \ref{th_Dirichlet}. 
\end{demo}

\begin{remark}\label{BC_robin}
In \eqref{eq_robin_2}, the boundary conditions satisfied by $v$ are not explicitly stated in terms of $v$ but 
rather in terms of $u=h_\sharp v$. Even if it not necessary for proving Theorem \ref{th_robin}, it could be 
interesting to compute them completely.
The second statement of Proposition \ref{prop_div} implies
$$|J|\big\la\,\nu_0\,\big|\, J^{-1} \cdot h^*(D\cdot\grad_x (h_\sharp v))+ i\big(J^{-1} \cdot h^*(\tilde A_h 
(h_\sharp v))\big\ra 
 = 0~~\text{ on }\partial\Omega_0~. $$
Now, $h_\sharp v= h_*\frac{1}{\sqrt{|J|}}v$. Thanks to Proposition \ref{prop_grad}, the fact that $|J|$ 
never vanishes yields
\begin{equation*}
\Big\la (J^{-1})^t \nu_0 \Big| (h^*D)\cdot(J^{-1})^t\cdot \sqrt{|J|} \grad_y \big(\frac v{\sqrt{J}}\big) + i 
(h^*\tilde A_h) v \Big\ra
=0~~\text{ on }\partial\Omega_0~. 
\end{equation*}
When $\tilde A_h$ corresponds to $A_h=-\frac 12 h_*(\partial_t h)$ and $D(t,x)=Id$, we obtain the boundary 
condition presented in \eqref{condi_neumann2}. 
\end{remark}

%%%%%%%%%%%%%%%%%%%%%%%%%%%%%%%%%%%%%%%%%%%%%%%%%%%%%%%%%%%%%%%%%%%%%%%%%%%%%%%%
%%%%%%%%%%%%%%%%%%%%%%%%%%%%%%%%%%%%%%%%%%%%%%%%%%%%%%%%%%%%%%%%%%%%%%%%%%%%%%%%

\section{Moser's trick}\label{section_Moser}

The aim of this section is to prove the following result.
\begin{theorem}\label{th_moser_1}
Let $k\geq 1$ and $r\geq 0$ with $k\geq r$ and let $\alpha\in (0,1)$. Let $\Omega_0$ be a bounded connected open 
$\Cc^{k+2,\alpha}$ domain of $\RR^N$. Let $I\subset\RR$ be an interval and let $f\in 
\Cc^r(I,\Cc^{k-1,\alpha}(\overline\Omega_0,\RR^*_+))$ be such that $\int_{\Omega_0} f(t,y)\dd y=\meas(\Omega_0)$ 
for all $t\in I$. Then, there exists a family $(\varphi(t))_{t\in I}$ of diffeomorphisms of $\overline \Omega_0$ 
of 
class $\Cc^r(I,\Cc^{k,\alpha}(\overline\Omega_0,\overline\Omega_0))$ satisfying
\begin{equation}\label{eq_th1}
\left\{\begin{array}{ll}
\det (D_y\varphi(t,y))=f(t,y),~&y\in\Omega_0\\
\varphi(t,y)=y,&y\in\partial\Omega_0.\end{array}\right.
\end{equation}
\end{theorem}

As explained in Section \ref{section_intro}, an interesting change of variables for a PDE with moving 
domains would be a family of diffeomorphisms having Jacobian with constant determinant. This is the goal of 
Theorem \ref{th_moser_2}, which is a direct consequence of Theorem \ref{th_moser_1}.

\vspace{3mm}

{ \noindent \emph{\textbf{Proof of Theorem \ref{th_moser_2}:}}}
Let Theorem \ref{th_moser_1} be valid. We set $\tilde h(t) = h(t)\circ \varphi^{-1}(t)$ and we compute 
\begin{align*}
\det\big(D_y\tilde h(t,y)\big)&=\det\big(D_y h(t,\varphi^{-1}(y))\big)\det\big(D_y (\varphi^{-1})(t,y) \big)\\
&=\det\big(D_y h(t,\varphi^{-1}(y))\big)\det\big((D_y \varphi)^{-1}(t,\varphi^{-1}(y)\big)\\
&=\left(\frac{\det(D_y h)}{\det(D_y \varphi)}\right)(t,\varphi^{-1}(y)).
\end{align*}
Since we want to obtain a spatially constant right-hand side, the dependence in $\varphi^{-1}(y)$ is not 
important. Thus,
$$\forall y\in\Omega_0~,~~\det(D_y\tilde h(t,y))=\frac{\meas(\Omega(t))} {\meas(\Omega_0)}$$
if and only if 
$$\forall y\in\Omega_0~,~~\det(D_y \varphi(t,y))= \frac {\meas(\Omega_0)}{\meas(\Omega(t))} \det(D_y h(t,y))~.$$
To obtain such a diffeomorphism $\varphi$, it remains to apply Theorem \ref{th_moser_1} with $f(t,y)$ being the 
above right-hand side. {\hfill$\square$\\}

\subsection{The classical results}\label{section_Moser2}
Theorem \ref{th_moser_1} is an example of a family of results aiming to find diffeomorphisms having prescribed 
determinant of the Jacobian. Such outcomes are often referred as ``Moser's trick'' because they were originated by 
the famous work of Moser \cite{Moser}. A lot of variants can be found in the literature, depending on the needs of 
the reader. We refer to \cite{Dacorogna} for a review on the subject. The following result comes from 
\cite{CDK,Dacorogna-Moser}, see also \cite{Dacorogna}.
\begin{theorem} \label{th_DM}
{\bf Dacorogna-Moser (1990)}\\
Let $k\geq 1$ and $\alpha\in (0,1)$. Let $\Omega_0$ be a bounded connected open $\Cc^{k+2,\alpha}$ domain. Then, 
the 
following statements are equivalent.
\begin{enumerate}
\item[(i)] The function $f\in\Cc^{k-1,\alpha}(\overline\Omega_0,\RR^*_+)$ satisfies $\int_{\Omega_0} f = 
\meas(\Omega_0)$.
\item[(ii)] There exists $\varphi\in\Diff^{k,\alpha}(\overline\Omega_0,\overline\Omega_0)$ satisfying
$$\left\{\begin{array}{ll}
\det (D_y\varphi(y))=f(y),~&y\in\Omega_0\\
\varphi(y)=y,&y\in\partial\Omega_0.\end{array}\right.
$$
\end{enumerate}
Furthermore, if $c>0$ is such that $\max(\|f\|_\infty,\|1/f\|_\infty)\leq c$ then there exists a constant 
$C=C(c,\alpha,k,\Omega_0)$ such that 
$$\|\varphi-\id\|_{\Cc^{k,\alpha}}~\leq~C\|f-1\|_{\Cc^{k-1,\alpha}}~.$$
\end{theorem}

A complete proof of theorem \ref{th_DM} can be found in \cite{Dacorogna}. There, a discussion concerning 
the optimality of 
the regularity of the diffeomorphisms is also provided. In particular, notice that the 
natural gain of regularity from $f\in\Cc^{k-1,\alpha}$ to $\varphi\in\Cc^{k,\alpha}$ was not present in the first 
work of Moser \cite{Moser} and cannot be obtained through the original method. Cases with other space regularity 
are 
studied in \cite{RY,RY2}.

For $t$ fixed, Theorem \ref{th_moser_1} corresponds to Theorem \ref{th_DM}. Our main goal is to extend the result 
to a 
time-dependent measure $f(t,y)$. In the original proof of \cite{Moser}, Moser uses a flow method and constructs 
$\varphi$ by a smooth deformation starting at $f(t=0)=\id$ and reaching $f(t=1)=f$. In this proof, the smooth 
deformation is a linear interpolation $f(t)=(1-t)\id+t f$ and one of the main steps consists in solving an ODE 
with 
a non-linearity as smooth as $\partial_t f(t)$. The linear interpolation is of course harmless in such a context. 
But when we 
consider another type of time-dependence, we need to have enough smoothness to solve the ODE. Typically, at least 
$\Cc^{1}-$smoothness in time is required to use the
arguments from the original paper \cite{Moser} of Moser. Thus, the original proof of \cite{Moser} does not provide 
an optimal regularity, neither in space or time. In particular, in Theorem \ref{th_moser_2}, this type of proof 
provides a diffeomorphism $\tilde h$ with one space regularity less than $h$. We refer to \cite{Hiba,BK,Geiges} 
for 
other time-dependent versions 
of Theorem \ref{th_DM}.

Our aim is to obtain a better regularity in space and time as well as to provide a complete proof of a 
time-dependent version of Moser's trick. To this end, we follow a method coming from 
\cite{Dacorogna,Dacorogna-Moser}, which is already known for obtaining the optimal space regularity. This method 
uses a fixed point argument, which can easily be parametrized with respect to time. Nevertheless, the fixed point 
argument only provides a local construction which is difficult to extend by gluing several similar construction 
(equations as \eqref{eq_th1} have an infinite number of solutions and the lack of uniqueness makes difficult to 
glue smoothly the different curves). That is the reason why we follow a similar strategy to the ones adopted in 
\cite{Dacorogna,Dacorogna-Moser}. First, at the cost of losing some regularity, we
prove the global result with the flow method. Second, we exploit a fixed 
point 
argument in order to ensure the result with respect to the optimal regularity.

\subsection{The flow method}
By using the flow method of the original work of Moser \cite{Moser}, we obtain the following version of Theorem 
\ref{th_moser_1}, where the statements on the regularities are weakened. 
\begin{prop}\label{prop1}
Let $k\geq 2$ and $r\geq 1$ with $k\geq r$. Let $\Omega_0$ be a bounded connected open $\Cc^{k,\alpha}$ domain of 
$\RR^N$ for some $\alpha\in (0,1)$. Let $I\subset\RR$ be an interval of times and let $f\in 
\Cc^r(I,\Cc^{k-1}(\overline\Omega_0,\RR^*_+))$ be such that $\int_{\Omega_0} f(t,y)\dd y=\meas(\Omega_0)$ for all 
$t\in I$. Then, there exists a family $(\varphi(t))_{t\in I}$ of diffeomorphisms of $\overline \Omega_0$ of class 
$\Cc^r(I,\Cc^{k-1}(\overline\Omega_0,\overline\Omega_0))$ satisfying \eqref{eq_th1}. 

Moreover, there exist continuous functions $M\mapsto C(M)$ and $M\mapsto \lambda(M)$ such that, if
$$\forall t\in I~,~~\|f(t,\cdot)\|_{\Cc^{k-1}}~\leq ~ M\text{ and }\min_{y\in\Omega_0} f(t,y) \geq \frac 1{M},$$
then 
$$\forall t\in I~,~~\|\varphi(t,\cdot)\|_{\Cc^{k-1}}\,+\,\|\varphi^{-1}(t,\cdot)\|_{\Cc^{k-1}} \,\leq\, C(M) 
e^{\lambda(M)|t|}.$$
\end{prop}
\begin{demo}
Let $t_0\in I$. Due to theorem \ref{th_DM}, there is $\varphi(t_0)$ with the required space regularity satisfying 
\eqref{eq_th1} at $t=t_0$. Setting $\varphi(t)=\tilde\varphi(t)\circ\varphi(t_0)$, we replace \eqref{eq_th1} by the 
condition
$$\det(D_y\tilde\varphi(t,y))=\frac{f(t,\varphi^{-1}(t_0,y))}{f(t_0,\varphi^{-1}(t_0,y))}~.$$
Thus, we may assume without loss of generality that $f(t_0,y)\equiv 1$. 

Let $L^{-1}$ be the right-inverse of the divergence introduced in Appendix \ref{section-inv-div}. Notice that 
$\partial_t f(t,\cdot)$ is of class $\Cc^{k-1}$ and hence of class $\Cc^{k-2,\alpha}$. Since 
$\int_{\Omega_0}\partial_t f=\partial_t \int_{\Omega_0} f=0$ we get that $f(t) \in Y^{k-2,\alpha}_m$ for all $t\in 
I$ and we can define 
\begin{equation}\label{eq_preuve_moser1}
U(t,\cdot)~=~-\frac 1{f(t,\cdot)} L^{-1}(\partial_t f(t,\cdot))~.
\end{equation}
We get that $U$ is well defined and it is of class $\Cc^{r-1}$ in time and $\Cc^{k-1}$ in space. For 
$x\in\overline\Omega_0$, we define $t\mapsto\psi(t,x)$ as the flow corresponding to the ODE
\begin{equation}\label{eq_preuve_moser2}
\psi(t_0,x)=x~~\text{ and }~~\partial_t \psi(t,x)=U(t,\psi(t,x))~~~~t\in I~.
\end{equation}
Notice that $\psi$ is locally well defined because $(t,y)\mapsto U(t,y)$ is at least lipschitzian in space and 
continuous in time. Moreover, the trajectories $t\in I\mapsto \psi(t,x)\in\overline\Omega_0$ are globally defined 
because $U(t,y)=0$ on the boundary of $\Omega_0$, providing a barrier of equilibrium points. The classical 
regularity results for ODEs show that $\psi$ is of class $\Cc^{r}$ with respect to time and $\Cc^{k-1}$ in space 
(see Proposition \ref{prop_regflow} in the appendix). Moreover, by reversing the time, the flow of a classical ODE 
as \eqref{eq_preuve_moser2} is invertible and $\psi(t,\cdot)$ is a diffeomorphism for all $t\in I$. We set 
$\varphi(t,\cdot)=\psi^{-1}(t,\cdot)$. 

Using Proposition \ref{prop_jacobi}, we compute for $t\in I$
\begin{align*}
\partial_t \big[\det(D_y\psi(t,y))\,f(t,\psi(t,y))\big]&= (\partial_t\det(D_y\psi(t,y)))\,f\circ\psi(t,y) + 
\det(D_y\psi(t,y)) \, (\partial_t f)\circ\psi(t,y)\\ &~~~~~~ + \det(D_y\psi(t,y)) ((D_y f)\circ\psi) \cdot 
(\partial_t \psi)(t,y)\\
&=\det(D_y\psi)\Big[ \Tr\big((D_y\psi)^{-1} \partial_t (D_y\psi)\big)f\circ\psi\\ &~~~~~~ +   
(\partial_t f)\circ\psi  +  \big((D_y f)\circ\psi\big) \cdot (\partial_t \psi) \Big](t,y)~.
\end{align*}
Since we have $D_y(\rho\circ\varphi)\circ\psi= (D_y\rho)(D_y\psi)^{-1}$ for any function $\rho$,
we use the trick
$$\Tr\big((D_y\psi)^{-1} \partial_t D_y(\psi)\big)=\Tr\big(D_y(\partial_t\psi)(D_y\psi)^{-1}\big) = 
\Tr\Big[D_y\big((\partial_t\psi)\circ\varphi\big)\Big]\circ\psi= \div_y 
\big((\partial_t\psi)\circ\varphi\big)\circ\psi.$$
Since $(\partial_t \psi)(t,\varphi(t,y))=U(t,\psi\circ\varphi(t,y))=U(t,y)$ due to \eqref{eq_preuve_moser2}, we 
obtain
\begin{align}\label{last}
\partial_t \Big[\det\big(D_y\psi(t,y)\big)\,f\big(t,\psi(t,y)\big)\Big]&= \det\big(D_y\psi(t,y)\big)\big[ 
\div_y(U)f + (\partial_t f) + (D_y f)U  \big]\circ\psi\\
&=\det\big(D_y\psi(t,y)\big)\big[ \div_y(Uf) + (\partial_t f) \big]\circ\psi~.
\end{align}
It remains to notice that $\div_y(Uf)=-\partial_t f$ by \eqref{eq_preuve_moser1}, that $f(t,y_0)=1$ and that 
$\det(D_y\psi(t_0,x))=\det(\id)=1$. Now,  from \eqref{last}, we have
$\det\big(D_y\psi(t,y)\big)\,f\big(t,\psi(t,y)\big)=\det\big(D_y\psi(t_0,y)\big)\,f\big(t,\psi(t_0,y)\big)=1$ and 
then
\begin{equation}\label{eq_preuve_moser3}
\forall t\in I,x\in\Omega_0~,~~  \det(D_y\psi(t,y))=\frac{1}{f(t,\psi(t,y))}~.
\end{equation}
We conclude that, as required,  
$$\forall t\in 
I,x\in\Omega_0~,~~\det(D\varphi(t,y))=\det((D\psi)^{-1}(t,\varphi(t,y)))=f(t,\psi\circ\varphi(t,y))=f(t,y)~. $$
Finally, it remains to notice that the bounds on $f$ yield bounds on $U$ due to \eqref{eq_preuve_moser1}. The 
classical bounds on the flow recalled in Proposition \ref{prop_regflow} in the appendix provide the claimed bounds 
on 
$\varphi$ and $\varphi^{-1}=\psi$. Indeed, they satisfy \eqref{eq_preuve_moser2} and the dual ODE where the 
time is reversed.  
\end{demo}

The trick of the above proof is the use of the flow of the ODE \eqref{eq_preuve_moser2}, which has of course a 
geometrical background. Actually, this idea can be extended to other differential forms than the volume form 
considered here, see \cite{Dacorogna} and \cite{Moser}. Also notice the explicit bounds stated in our version, 
which is not usual. We will need them for the proof of Theorem \ref{th_moser_1} because we slightly adapt the 
strategy of \cite{Dacorogna-Moser}.

\subsection{The fixed point method}
Proposition \ref{prop2} below is an improvement of Proposition \ref{prop1} regarding regularity in both time and 
space, but it only deals with local perturbations of $f(t,y)\equiv 1$. It is proved following the fixed point 
method of \cite{Dacorogna,Dacorogna-Moser}. 
\begin{prop}\label{prop2}
Let $k\geq 1$ and $r\geq 0$ with $k\geq r$ and let $\alpha\in (0,1)$. Let $\Omega_0$ be a bounded connected open 
$\Cc^{k+1,\alpha}$ domain of $\RR^N$. There exists $\varepsilon>0$ such that, for any interval $I\subset\RR$ and 
any $f\in \Cc^r(I,\Cc^{k-1,\alpha}(\overline\Omega_0,\RR^*_+))$ being such that 
$$\forall t\in I~,~~\int_{\Omega_0} f(t,y)\dd y=\meas(\Omega_0)~~\text{ and 
}~~\|f(t,\cdot)-1\|_{\Cc^{k-1,\alpha}}\leq \varepsilon~.$$
Then, there exists a family $(\varphi(t))_{t\in I}$ of diffeomorphisms of $\overline \Omega_0$ of class 
$\Cc^r(I,\Cc^{k,\alpha}(\overline\Omega_0,\overline\Omega_0))$ satisfying \eqref{eq_th1}. Moreover, there exists a 
constant $K$ independent of $f$ such that
$$\sup_{t\in I} \|\varphi(t)\|_{\Cc^{k,\alpha}} + \sup_{t\in I}\|\varphi^{-1}(t)\|_{\Cc^{k,\alpha}} 
~\leq~K\sup_{t\in I} \|f(t,\cdot)-1\|_{\Cc^{k-1,\alpha}}~.$$
\end{prop}

\begin{demo}
We follow Section 5.6.2 of \cite{Dacorogna}. We reproduce the proof for sake of completeness and to explain why 
we can add for free the time-dependence. We set 
$$Q~:~M\in\Mc_{d}(\RR)~\longmapsto~\det(\id+M)-1-\Tr(M)\in\RR~.$$
Notice that $Q$ is the sum of monomials of degree between $2$ and $d$ with respect to the coefficients because $1$ 
and $\Tr(M)$ are the first order terms of the development of $\det(\id+M)$. Using the fact that $\Cc^{k,\alpha}$ 
is 
a Banach algebra (see for example \cite{CDK}), there exists a constant $K_2$ such that, for any functions 
$u,v\in\Cc^{k,\alpha}(\overline\Omega,\RR^N)$
\begin{align}
\|Q(D_yu)-Q(D_yv)\|_{\Cc^{k-1,\alpha}}\leq K_1(1+&\|u\|_{\Cc^{k,\alpha}}^{d-2}+\|v\|_{\Cc^{k,\alpha}}^{d-2})
\nonumber\\ &\max(\|u\|_{\Cc^{k,\alpha}},\|v\|_{\Cc^{k,\alpha}})\|u-v\|_{\Cc^{k,\alpha}}~.\label{eq_proof_th2_3}
\end{align}
We seek for a function $\varphi(t,y)$ solving \eqref{eq_th1} for times $t$ close to $t_0$ by setting 
$$\varphi(t,y)= \id + \eta(t,y)~.$$
Since we would like that $\det(\id+D_y\eta)=f$ and as $\div\eta = \Tr(D_y\eta)$, the identity \eqref{eq_th1} is 
equivalent to 
\begin{equation}\label{eq_proof_th2_1}
\left\{\begin{array}{ll}
\div \eta(y,t)=f(t,y)-1-Q(D_y\eta(t,y))~&y\in\Omega_0\\
\eta(t,y)=0&y\in\partial\Omega_0\end{array}\right.
\end{equation}
Let $X^{k,\alpha}_0$ and $Y^{k-1,\alpha}_m$ the spaces introduced in Appendix \ref{section-inv-div}. First notice 
that, by assumption, 
$$\int_{\Omega_0} f(t,y)\dd y = \meas(\Omega_0)= \int_{\Omega_0} 1 \dd y$$
and thus $(f(t)-1)$ belongs to $Y^{k-1,\alpha}_m$. 
Since $f(t,\cdot)$ is close to $1$, we seek for $\eta$ small in $\Cc^r(I,\Cc^{k,\alpha}(\overline\Omega_0,\RR^N))$. 
By Corollary \ref{coro_MO_2}, there exists $R>0$ small such that if $\|\eta(t,\cdot)\|_{\Cc^{k,\alpha}}\leq R$ for 
all $t\in I$, then $\varphi(t,y)= \id + \eta(t,y)$ is a diffeomorphism from $\Omega_0$ onto itself. If this is the 
case, we have that 
\begin{align*}
\int_{\Omega_0} Q(D_y\eta(t,y)) \dd y &= \int_{\Omega_0} \big(\det(\id+D_y\eta)(t,y)-1+\Tr(D_y\eta(t,y)) \big)\dd 
y 
\\
&=\int_{\Omega_0}\det(D_y\varphi(t,y))\dd y -\meas(\Omega_0) + \int_{\Omega_0}\div(\eta(t,y))\dd y\\
&=\int_{\Omega_0} \Un(\varphi(y))\dd y -\meas(\Omega_0) + \int_{\partial\Omega_0}\eta(y)\dd y\\
&= \meas(\Omega_0) - \meas(\Omega_0) + 0 = 0~.
\end{align*}
This means that we can look for $\eta$ small as a solution in $\Cc^r(I,\Cc^{k,\alpha}(\overline\Omega_0,\RR^N))$ 
of 
\begin{equation}\label{eq_proof_th2_2}
\eta(t)\,=\,L^{-1}\big(f(t)-1-Q(D_y\eta(t))\big)
\end{equation}
since we expect $Q(D_y\eta(t))$ to belong to $Y^{k-1,\alpha}_m$. We construct $\eta(t)$ by a fixed point argument 
by applying Theorem \ref{th_fixed_point} to the function 
$$\Phi~:~(\eta,t)\in \Bc_R\times I~\longmapsto~L^{-1}\big(f(t)-1-Q(D_y\eta)\big)\in X^{k,\alpha}_0$$
where $L^{-1}$ is the right-inverse of the divergence (see Appendix \ref{section-inv-div}) and 
$\Bc_R=\{\eta\in X^{k,\alpha}_0, \|\eta\|_{X^{k,\alpha}}\leq R\}$ with $R$ as above. As a consequence,  
$\id+\eta$ is a diffeomorphism and the above computations are valid. Due to \eqref{eq_proof_th2_3}, we have 
$$\|\Phi(\eta_1,t)-\Phi(\eta_2,t)\|_{\Cc^{k,\alpha}}~\leq~ R K_1K_2 
(1+2R^{d-2})\|\eta_1-\eta_2\|_{\Cc^{k,\alpha}}$$
where $K_2=|||L^{-1}|||$. Using that $Q(0)=0$, we also have
$$\|\Phi(\eta,t)\|_{\Cc^{k,\alpha}}~\leq~K_2 \big( \|f(t)-1\|_{\Cc^{k-1,\alpha}} + K_1(R^2+R^{d}) \big)~.$$
We choose $R\in(0,1/2]$ small enough such that $\id+\eta$ is a diffeomorphism and  
$$RK_1K_2(1+2R^{d-2})~\leq~ \frac 12$$
which also implies that $K_1K_2(R^2+R^{d})\leq R/2$. Then, we choose $\varepsilon=R/2K_2$ and assume that 
\begin{equation}\label{eq_proof_th2_6}
\forall t\in I~,~~\|f(t)-1\|_{\Cc^{k-1,\alpha}}~\leq~\varepsilon=\frac{R}{2K_1}~.
\end{equation}
By construction, $\Phi$ is $1/2-$lipschitzian from $\Bc_R$ into $\Bc_R$. The classical fixed point theorem for 
contraction maps shows the existence of $\eta(t)$ solving \eqref{eq_proof_th2_2} for each $t$. The regularity of 
$\eta$ with respect to the time $t$ is directly given by Theorem \ref{th_fixed_point}. By construction 
$\eta(t)=L^{-1}(f(t)-1-Q(D_y\eta(t))$, and since $LL^{-1}=\id$, this implies that $\div 
\eta(t)=f(t)-1-Q(D_y\eta(t))$ and thus $f(t)=\det(\id+\eta(t))=\det(\varphi(t))$. Also remember that, by 
construction,  $\varphi$ is $\Cc^k$ close enough to the identity in $\Omega_0$ and is the identity on 
$\partial\Omega_0$ so that $\varphi(t)$ is a $\Cc^k-$diffeomorphism (see Appendix \ref{appendix_MO} and the 
topological arguments of \cite{Meisters-Olech}).
\end{demo}

\subsection{Proof of Theorem \ref{th_moser_1}}
As in the work of Dacorogna and Moser (see \cite{Dacorogna} and \cite{Dacorogna-Moser}), we obtain Theorem 
\ref{th_moser_1} by combining  the propositions \ref{prop1} and \ref{prop2}. In the original method, the authors 
used first the fixed point method and then the flow method. However, this approach does not provide the optimal 
time regularity. Therefore, we couple the two methods in the other 
sense. In this case, we need the $\Cc^{k-1}$-bounds on the diffeomorphism provided by Proposition \ref{prop1}. 
This is the reason why we made them explicit, which is not common for this type of results.

Let $f\in\Cc^r(I,\Cc^{k-1,\alpha}(\overline\Omega_0,\RR^*_+))$ with $\int_\Omega f(t,y)\dd y=\meas(\Omega_0)$ for 
all $t\in I$. Let $\varepsilon>0$. We consider a regularization $f_1$ of $f$, which is of class 
$\Cc^{r'}(I,\Cc^{k'-1}(\overline\Omega_0,\RR^*_+))$ with $r'=\max(1,r)$ and $k'=k+2$, satisfying
\begin{equation}\label{eq_proof_th1}
\forall t\in I~,~~\left\|\frac{f(t,\cdot)}{f_1(t,\cdot)}-1\right\|_{\Cc^{k-1,\alpha}} \leq \varepsilon(t)~.
\end{equation}
Assume that we first apply the flow method of Proposition \ref{prop1} to obtain a smooth family of diffeomorphisms 
$\varphi_1\in \Cc^{r'}(I,\Cc^{k'-1}(\overline\Omega_0,\overline\Omega_0))\subset 
\Cc^{r}(I,\Cc^{k+1}\overline\Omega_0,\overline\Omega_0))$ such that 
$$\det(D_y\varphi_1(t,y))=f_1(t,y)~~\text{ and }~~\varphi_{1|\partial\Omega_0}(t,\cdot)=\id_{|\partial\Omega_0}~.$$
Now, we seek for $\varphi$ of the form $\varphi=\varphi_2\circ\varphi_1$ satisfying the statement of Theorem 
\ref{th_moser_1}. Thus, we need to find $\varphi_2\in \Cc^{r}(I,\Cc^{k,\alpha}(\overline\Omega_0,\RR^*_+))$ 
such that 
$$\det(D_y\varphi_2)(t,\varphi_1(t,y)){\det(D_y\varphi_1(t,y))}={f(t,y)}$$
{\it i.e.}
\begin{equation}\label{eq_proof_th1_2}
\det(D_y\varphi_2)(t,y)=\frac{f(t,\varphi_1^{-1}(t,y))}{f_1(t,\varphi_1^{-1}(t,y))}~.
\end{equation}
We would like to apply the fixed point method of Proposition \ref{prop2} by choosing $\varepsilon(t)$ small enough 
such that
\begin{equation}\label{eq_proof_th1_3}
\forall t\in I~,~~\left\|\frac{f(t,\varphi_1^{-1}(t,y))}{f_1(t,\varphi_1^{-1}(t,y))}-1\right\|_{\Cc^{k-1,\alpha}} 
\leq \varepsilon
\end{equation}
with $\varepsilon$ as in Proposition \ref{prop2}. However, we notice that in \eqref{eq_proof_th1_3} appears the 
composition with $\varphi_1^{-1}$. Fix $\varepsilon$ as in Proposition \ref{prop2} and smaller than $1/2$. When we 
consider Proposition \ref{prop1} with $M=2$, we have the upper bound for $C(2)e^{\lambda(2)t}$ for the 
derivatives of any $\varphi_1$ constructed as above when $\varepsilon(t)\leq \varepsilon$ in \eqref{eq_proof_th1}. 
In particular, these bounds indicate how much the composition by $\varphi^{-1}$ increases the 
$\Cc^{k-1,\alpha}-$norm. 
Having this in mind, we may choose $\varepsilon(t)$ small (and exponentially decreasing) such that, for any $f_1$ 
satisfying \eqref{eq_proof_th1}, the associated diffeomorphisms $\varphi_1$ are such that \eqref{eq_proof_th1_3} 
holds. 

We construct the diffeomorphism as follows. We choose a regularization $f_1$ of the function $f$, which is of 
class $\Cc^{r'}(I,\Cc^{k'-1}(\overline\Omega_0,\RR^*_+))$ with $r'=\max(1,r)$ and $k'=k+2$, satisfying 
\eqref{eq_proof_th1} for the suitable $\varepsilon(t)$. We will also need that $\int_{\Omega_0}(f/f_1)(t,y)\dd 
y=\meas(\Omega_0)$, which is provided by first choosing $\tilde f$ satisfying \eqref{eq_proof_th1} with a smaller 
$\varepsilon(t)$ and then set $f_1=\tilde f\cdot\int(f/\tilde f)/\meas(\Omega_0)$. By Proposition \ref{prop1}, 
there exists a family of diffeomorphisms $t\mapsto\varphi_1(t,\cdot)$ satisfying \eqref{eq_proof_th1} and of class 
$\Cc^{r'}(I,\Cc^{k'-1}(\overline\Omega_0,\overline\Omega_0))\subset\Cc^{r}(I,\Cc^{k+1}(\overline\Omega_0,
\overline\Omega_0))$. 

Now we consider $f_2=f/f_1$ which satisfies by construction \eqref{eq_proof_th1_3} with $\varepsilon$ as in 
Proposition \ref{prop2}. Notice that $f_2$ is of class $\Cc^r(I,\Cc^{k-1,\alpha}(\overline\Omega_0,\RR^*_+))$ and 
\begin{align*}
\int_{\Omega_0} f_2(t,y)\dd y&= \int_{\Omega_0} 
\frac{f(t,\varphi_1^{-1}(t,y))}{\det(D_y\varphi_1(t,\varphi_1^{-1}(t,y)))}\dd y\\
&=\int_{\Omega_0}{f(t,\varphi_1^{-1}(t,y))}{\det(D_y(\varphi_1^{-1}(t,y)))}\dd y\\
&=\int_{\Omega_0}f(t,y)\dd y =\meas (\Omega_0)~.
\end{align*}
Thus, we can apply Proposition \ref{prop2} to obtain a family of diffeomorphisms $\varphi_2$ of class 
$\Cc^r(I,\Cc^{k,\alpha}(\overline\Omega_0,\overline\Omega_0))$ such that
$$\det(D_y\varphi_2(t,y))=f_2(t,y)~~\text{ and }~~\varphi_{2|\partial\Omega_0}(t,\cdot)=\id_{|\partial\Omega_0}~.$$
By construction, $\varphi=\varphi_2\circ\varphi_1$ satisfies the conclusion of Theorem \ref{th_moser_1}.

%%%%%%%%%%%%%%%%%%%%%%%%%%%%%%%%%%%%%%%%%%%%%%%%%%%%%%%%%%%%%%%%%%%%%%%%%%%%%%%%
%%%%%%%%%%%%%%%%%%%%%%%%%%%%%%%%%%%%%%%%%%%%%%%%%%%%%%%%%%%%%%%%%%%%%%%%%%%%%%%%

\section{Some applications of our results}\label{section_appli}
\subsection{Adiabatic dynamics for quantum states on moving domains}
In this subsection, we show how to ensure an adiabatic result for 
the Schr\"odinger equation on a moving domain as \eqref{SE_Dirichlet_2}. 
We consider the framework of Corollary \ref{coro_adiabatic}. 
We denote by $H(\tau)=-\Delta$ the Dirichlet Laplacian in $L^2(\Omega(\tau),\CC)$, {\it i.e.} 
$D(H(\tau))=H^2(\Omega(\tau))\cap H^1_0(\Omega(\tau)) $ and $\tau\in[0,1]\mapsto\lambda(\tau)$, a continuous curve 
such that $\lambda(\tau)$ for every $\tau\in [0,1]$ is in the discrete spectrum of $H(\tau)$. We also assume that 
$\lambda(\tau)$ is a simple isolated eigenvalue for every $\tau\in[0,1]$, associated with the spectral projectors 
$P(\tau)$. 

We consider on the time interval $[0,1/\epsilon]$ the equation \eqref{eq_dirichlet_3} in $L^2(\Omega_0,\CC)$ when 
$H(t)$ is a Dirichlet Laplacian. Fixed $\epsilon>0,$ we 
substitute $t$ by $\tau=\epsilon t$ and set $\tilde v_\epsilon(\tau)=v(\tau/\epsilon)=h^\sharp 
u_\epsilon(\tau/\epsilon)$ to obtain 
\begin{equation}\begin{split}\begin{cases}\label{SE_adiabatic2}
i\epsilon\ddd_\tau \tilde v_\epsilon(\tau) =\Big(h^\sharp 
H(\tau) h_\sharp+\epsilon \Hc(\tau)\Big) \tilde v_\epsilon(\tau), &\ \ \ \ \  \ \ \ \tau\in [0,1]\\
\tilde v_\epsilon(\tau)_{|\partial\Omega_0}\equiv 0,& \\
\tilde v_\epsilon(\tau=0)=h^\sharp u_0 &
\end{cases}\end{split}\end{equation} 
where $$\Hc(\tau)=- h^\sharp\Big[i \div_x \big(A_h(\tau,x) \cdot \big) ~+~ i \big\la A_h(\tau,x)\big| 
\grad_x \cdot \big\ra\Big] h_\sharp~,\ \ \ \ \ \  A_h(\tau,x):=-\frac 12 (h_*\partial_\tau h)(\tau,x).$$  
The problem \eqref{SE_adiabatic2} generates a unitary flow thanks to Section \ref{section_cauchy}. Even though it 
is well known that the classical adiabatic theorem is valid for the dynamics 
$i\epsilon\ddd_\tau u_\epsilon(\tau) =h^\sharp 
H(\tau) h_\sharp u_\epsilon(t)$, (see \cite[Chapter\ 4]{Bornemann} or 
\cite{Teufel}), we may wonder if it is the same for the equation \eqref{SE_adiabatic2} because 
$$\tilde H_\epsilon(\tau)=h^\sharp H(\tau) h_\sharp + \epsilon \Hc(\tau) $$ 
depends on $\epsilon$ and no spectral assumptions have been made on this family. First, we notice that, by 
conjugation, $\lambda(\tau)$ also belongs to the discrete spectrum of the operator $h^\sharp H(\tau) h_\sharp$ in 
$L^2(\Omega_0,\CC)$ and is associated with the spectral projection $(h^\sharp P h_\sharp)(\tau)$. 
Then, for each $\tau$, $\tilde H_\epsilon(\tau)$ is a small relatively compact self-adjoint perturbation of 
$h^\sharp H(\tau) h_\sharp$. Thus, for all $\epsilon>0$ small enough, there exists a curve 
$\tilde\lambda_\epsilon(\tau)$ of simple isolated eigenvalues of $\tilde H_\epsilon(\tau)$, associated with spectral 
projectors $\tilde P_\epsilon(\tau)$, such that $\tilde \lambda_\epsilon$ and $\tilde P_\epsilon$ converge 
uniformly when $\epsilon\rightarrow 0$ to $\lambda$ and $h^\sharp P h_\sharp$ respectively 
 (see \cite{Kato} for further details). In this framework, even if $\tilde H_\epsilon(\tau)$ depends on 
$\epsilon$, 
it is known that the classical adiabatic arguments can still be applied (see for example Nenciu \cite[Remarks;\ 
p.\ 16;\ (4)]{Nen}, Teufel \cite[Theorem\ 4.15]{Teufel} or the works \cite{ASY, Ga, JoPf}).
 Thus, we obtain the following convergence for the solution of \eqref{SE_adiabatic2}
$$\la \tilde P_\epsilon(1)\tilde v_\epsilon(1)|\tilde v_\epsilon(1)\ra ~~\xrightarrow[~~\epsilon\longrightarrow 
0~~]{} ~~\la \tilde P_0(0)v(0)|v(0)\ra~=~\la (h^\sharp P h_\sharp)(0).h^\sharp u_0|h^\sharp u_0\ra~=~ \la P(0) u_0| 
u_0\ra~.$$
Finally, we notice that, since $\tilde P_\epsilon(1)$ converges to $h^\sharp P(1) h_\sharp$ when $\epsilon$ goes to 
zero, 
$$\la \tilde P_\epsilon(1)\tilde v_\epsilon(1)|\tilde v_\epsilon(1)\ra~=~\la (h^\sharp P(1) h_\sharp) \tilde 
v_\epsilon(1) |\tilde  v_\epsilon(1)\ra ~+~o(1)~=~\la P(1) u_\epsilon(1/\epsilon)|u_\epsilon(1/\epsilon)\ra 
~+~o(1)~.$$
This concludes the proof of Corollary \ref{coro_adiabatic}.

\subsection{Explicit examples of time-varying domains}
\needspace{5mm}
\noindent{\bf \underline{Translation of a }p\underline{otential well}}\\[1mm]
Let us consider any domain $\Omega_0\subset\RR^N$ and any smooth family of vectors $D(t)\in \Cc^2([0,T],\RR^N)$. 
The family of translated domains is $\Omega(t)=h(t,\Omega_0)$ where $h(t,y)=y+D(t)$. By explicit computations, we 
obtain that $h^\sharp \Delta_x h_\sharp = \Delta_y$ and $(h_*\ddd_th)(t,x)=D'(t)$. 
Since $|J|$ does not depend on $y$, we do not need Moser's trick to get \eqref{SE_moser}. We can apply the 
gauge transformation since $\phi(t,x)=\frac{1}{2} \la D(t)|x\ra$ satisfies $h_*\ddd_th=2\grad_x\phi$. Then, 
$w=h^\sharp e^{-i\phi}u$ satisfies Equation \eqref{SE_Gauge}, which becomes in our framework
\begin{equation}\label{SE_Gauge_example1}
i\partial_t w~=~-\Delta_y w ~+~ \frac{1}{4}\Big( 2\la D''(t)|(y+D(t))\ra + |D(t)|^2 \Big) w.
\end{equation}
In this very particular case, we can further simplify the expression thanks to an interesting fact. Two 
terms of the electric potential in \eqref{SE_Gauge_example1} do not depend on the space variable. We can thus 
apply 
another transformation to the system by adding a phase which is an antiderivative of $2\la D''(t)|D(t)\ra + 
|D(t)|^2$. For example, we consider 
$$\tilde w~=~e^{\frac i4( 2\la D'(t)|D(t)\ra-\int_0^t |D(s)|^2\dd s)}w$$ 
where $w$ satisfies Equation \eqref{SE_Gauge_example1}. Then, $\tilde w$ is solution of the equation
\begin{equation}\label{SE_Gauge_example1bis}
i\partial_t \tilde w~=~-\Delta_y \tilde w ~+~ \frac{1}{2}\la D''(t)|y \ra \tilde w.
\end{equation}
These explicit computations are not new and appear for example in \cite{BeauchardCoron} for the one dimensional 
case of a moving interval.

\vspace{5mm}

\needspace{5mm}
\noindent{\bf \underline{Rotatin}g\underline{ domains }}\\[1mm]
Let us consider a family of rotating domains in $\RR^2$. Clearly, the same results can be extended by considering 
rotations in $\RR^N$ with $N\geq 3$.
Let $\Omega_0\subset\RR^2$ and let $\Omega(t)=h(t,\Omega_0)$ with $$h(t):~y=(y_1,y_2) \in \Omega_0~\longmapsto~ 
\big(\cos(\omega t )y_1-\sin(\omega t)y_2,\cos(\omega t )y_1+\sin(\omega t)y_2\big)$$
and $\omega\in\RR$. Using the classical notation $(y_1,y_2)^\perp=(-y_2,y_1)$, it is straightforward to check that 
$|J|=1$, $J^{-1}(J^{-1})^t=Id$, 
$$(h_*\partial_t h)(t,x)=\omega x^\perp~~\text{ and }~~J^{-1}\ddd_t h(t,y)= \omega y^\perp~.$$
We obtain by direct computation or by the first line of \eqref{SE_moser} that $v=h^*u$ satisfies the following 
Schr\"odinger equation in the rotating frame
\begin{equation}\label{SE_Gauge_example1ter}
i\partial_t v~=~-\Delta_y v ~+~ i\omega\la y^\perp,\grad_y v\ra~.
\end{equation}
This is an obvious and well-known computation (used for studying quantum systems in rotating potentials frames). 
The general Hamiltonian structure highlighted in this paper simply writes here as 
 
$$i\partial_t v~=~-\Big(\grad_y -i\frac \omega 2 y^\perp\Big)^2v ~-~ \frac {\omega^2}4 |y|^2 v ~. $$
which is given by the second line of \eqref{SE_moser} or obtained directly from \eqref{SE_Gauge_example1ter} by 
using the fact that $y^\perp$ is divergence free. We recover a repulsive potential $\omega^2|y|^2/4$ corresponding 
to the centrifugal force.

\vspace{5mm}

\needspace{5mm}
\noindent{\bf \underline{Movin}g\underline{ domains with dia}g\underline{onal diffeomor}p\underline{hisms}}\\[1mm]
Let $\Omega_0\subset\RR^N$ and let $\Omega(t)=h(t,\Omega_0)$ with $$h(t):~y=(y_i)_{i=1\ldots N} \in 
\Omega_0~\longmapsto~ 
\big(f_i(t)\,y_i\big)_{i=1\ldots N}$$
and $f_i \in\Cc^2([0,T],\RR)$. As above, we obtain 
\begin{equation}\label{eq_example2}
h^\sharp \Delta_x h_\sharp = \sum_{i=1}^N \frac{1}{f_i(t)^2}\ddd_{y_iy_i}^2,\ \ \ \ \ (h_*\ddd_th)(t,x)=\left( 
\frac{f_i'(t)}{f_i(t)}\,x_i\right)_{i=1\ldots N}.
\end{equation}
We apply again the gauge 
transformation and then,
\begin{equation}\label{eq_example}\phi(t,x)=\frac{1}{4}\sum_{j=1}^N\left(\frac{f_i'(t)}{f_i(t)}\,x_i^2\right)\ \ \ 
\ \ \text{satisfies}\ \ \ \ \ h_*\ddd_th=2\grad_x\phi.\end{equation} Finally, $w=h^\sharp e^{-i\phi}u$ satisfies 
the equation
\begin{equation}\label{SE_Gauge_example}
i\partial_t w~=~-\sum_{i=1}^N \frac{1}{f_i(t)^2}\ddd_{y_iy_i}^2 w ~+~ \frac{1}{4}\left(\sum_{i=1}^N f''_i(t) 
f_i(t) y_i^2 \right) w.
\end{equation}
In the homothetical case where $f_1=f_2=\ldots =: f(t)$ 
(see for instance 
\cite{BMT,Beauchard,BeauchardTeismann,Knobloch-Krechetnikov,Makowski-Peplowski,Moyano,Pinder,Rouchon}). Equation 
\eqref{SE_Gauge_example} becomes
\begin{equation*}
i\partial_t w(t)=-\frac{1}{f(t)^2}\Delta 
w(t)+ \frac{1}{4}f''(t) f(t) |y|^2 w(t),~~~t\in I.
\end{equation*}
In this case, it is usual to make a further simplification to eliminate the time-dependence of the main operator by 
changing the time variable for 
$$\tau=\int_0^t \frac{1}{f(s)^2}ds$$ 
and introducing the implicitly defined function 
$$U(\tau)=\frac{f'(t)f(t)}{4}.$$ 
We obtain
\begin{equation}\label{SE_Gauge_example_pulsating_2}
i\partial_\tau w(\tau)=-\Delta w(\tau)+\Big(U'(\tau)-4U(\tau)\Big)|y|^2 
w(\tau),~~~~~~~~~\tau\in\Big[\,0\,,\,\int_0^T 1/{f(s)^2}ds\Big].
\end{equation}
In this simple case, we see that the general framework of this paper coincides with the previous computations 
introduced in dimension $d=1$ by Beauchard in \cite{Beauchard}. Indeed, the transformations
described in \cite[Section\ 1.3]{Beauchard} corresponds to application $u\longmapsto w=h^\sharp e^{-i\phi}u$ as 
the multiplication for the square root of the exponential \cite[relation\ (1.4)]{Beauchard} is nothing else than 
the multiplication for the square root of the Jacobian appearing in the definition of $h^\sharp$. Our paper put the 
change of variable of \cite{Beauchard} in a more general geometric framework.

A similar expression to \eqref{SE_Gauge_example_pulsating_2} is also obtained by Moyano in \cite{Moyano} for the 
case of the two-dimensional disk, nevertheless the transformations adopted in \cite{Moyano} are different 
from the ones considered in our work. In particular, they are not unitary 
with respect to the classical $L^2$-norm.

\vspace{3mm}

For a second application of this simple case, we consider the case of a family of cylinders 
$$\Omega(t)=\Big\{(x_1,x_2,x_3)\in\RR_+^3\ \:\ x_1\in(0,\ell(t)),\ x_2^2+x_3^2<r^2,\ \Big\}$$
for $\ell$ a $\Cc^2-$varying length. 
We would like to consider the Schr\"odinger equation $i\partial_t u =-\Delta_x u$ in $\Omega(t)$ with boundary 
conditions of the Neumann type. 
This example corresponds to the situation of Figure \ref{fig-cylindre} in Section \ref{section_intro}. As shown in 
this paper, the conditions at the boundaries cannot be pure homogeneous Neumann ones everywhere if $\ell(t)$ is not 
constant. Theorem \ref{th_main_Robin} shows us the correct ones. We can choose 
$$h(t,\cdot):\big(y_1,y_2,y_3\big)\in \Omega_0\longmapsto 
\big(\ell(t)y_1,y_2,y_3\big)\in\Omega(t)$$
with $\Omega_0$ the cylinder of length $1$. Then, we get by \eqref{eq_example2}
the term $h_*\partial_t h$ and we can compute that the suitable boundary conditions are
\begin{equation}\label{eq_BC_cylindre}
\partial_\nu u-\frac i2 \ell'(t) u=0 ~~~ \text{ at }x\in\partial\Omega(t)\text{ with }x_1=\ell(t)
\end{equation}
and $\partial_\nu u=0$ on the other parts of the boundary (see Figure \ref{fig-cylindre}). It is important to 
notice that, contrary to the above changes of variables, this computation is independent of the choice of $h(t)$. 
Equations as  \eqref{SE_Gauge_example_pulsating_2} can be seen as auxiliary equations and they depend on several 
choices, while \eqref{eq_BC_cylindre} is stated for the original variable $u$ and have physical 
meaning.

%%%%%%%%%%%%%%%%%%%%%%%%%%%%%%%%%%%%%%%%%%%%%%%%%%%%%%%%%%%%%%%%%%%%%%%%%%%%%%%%
%%%%%%%%%%%%%%%%%%%%%%%%%%%%%%%%%%%%%%%%%%%%%%%%%%%%%%%%%%%%%%%%%%%%%%%%%%%%%%%%
%%%%%%%%%%%%%%%%%%%%%%%%%%%%%%%%%%%%%%%%%%%%%%%%%%%%%%%%%%%%%%%%%%%%%%%%%%%%%%%%
%%%%%%%%%%%%%%%%%%%%%%%%%%%%%%%%%%%%%%%%%%%%%%%%%%%%%%%%%%%%%%%%%%%%%%%%%%%%%%%%

\appendix

\section{Appendix}

\subsection{Unitary semigroups}\label{section-kisynski}
Defining solutions of an evolution equation with a time-dependent family of operators is nowadays a classical 
result (see \cite{Tanabe}). In the present article, we are interested in the Hamiltonian structure and we use the 
following result of \cite{Kisynski} (see also \cite{Teufel}).
\begin{theorem}\label{th_Kisynski}
{\bf (Kisy\'nski, 1963)}\\
Let $\Xc$ be a Hilbert space and let $(\Hc(t))_{t\in [0,T]}$ be a family of self-adjoint positive operators 
on $\Xc$ such that $\Xc^{1/2}=D(\Hc(t)^{1/2})$ is independent of time $t$. Also set 
$\Xc^{-1/2}=D(\Hc(t)^{-1/2})=(\Xc^{1/2})^*$ and assume that $\Hc(t):\Xc^{1/2}\rightarrow 
\Xc^{-1/2}$ is of class $\Cc^1$ with respect to $t\in [0,T]$. In other words, we assume that the sesquilinear form 
$\phi_t(u,v)=\la \Hc(t)u|v\ra$ associated with $\Hc(t)$ has a domain $\Xc^{1/2}$ independent of the time and is of 
class $\Cc^1$ with respect to $t$. Also assume that
there exist $\gamma>0$ and $\kappa\in\RR$ such that,
\begin{equation}\label{hyp-Kisynski}
\forall t\in [0,T]~,~\forall u\in \Xc^{1/2}~,~~ \phi_t(u,u)=\langle 
\Hc(t)u|u\rangle_{\Xc} \geq \gamma \|u\|^2_{\Xc^{1/2}} - \kappa \|u\|^2_{\Xc}~.
\end{equation}
Then, for any $u_0\in \Xc^{1/2}$, there is a unique solution $u$ belonging to $\Cc^0([0,T], \Xc^{1/2})\cap 
\Cc^1([0,T],\Xc^{-1/2})$ of the equation 
\begin{equation}\label{eq-Kisynski}
i\partial_t u(t)=\Hc(t)u(t)~~~~u(0)=u_0~.
\end{equation}
Moreover, $\|u_0\|_{\Xc}=\|u(t)\|_{\Xc}$ for all $t\in [0,t]$ and we may extend by 
density the flow of \eqref{eq-Kisynski} on $\Xc$ as a unitary flow $U(t,s)$ such 
that $U(t,s)u(s)=u(t)$ for all solutions $u$ of \eqref{eq-Kisynski}.

If in addition $\Hc(t):\Xc^{1/2}\rightarrow \Xc^{-1/2}$ is of class $\Cc^2$ and $u_0\in D(\Hc(0))$, then $u(t)$ 
belongs to $D(\Hc(t))$ for all $t\in [0,T]$ and $u$ is of class $\Cc^1([0,T],\Xc)$.
\end{theorem}

\subsection{The derivative of the determinant}
We recall the following standard result.
\begin{prop}\label{prop_jacobi}
Let $I\subset\RR$ be an interval of times and let $N\geq 1$. Let $A(t)$ be a family of $N\times N$ complex matrices 
which is differentiable with respect to the parameter $t\in I$. If $A(t)$ is invertible for every $t\in I$, then
$$\ddd_t \det(A(t))=\det(A(t))\Tr\big(A(t)^{-1}\ddd_t A(t)\big).$$
More generally, we have $\partial_t \det(A(t))=\Tr(\text{\rm com}(A(t))^t \partial_tA(t))$ where $\text{\rm 
com}(A)$ is the comatrix of $A$.
\end{prop}
\begin{demo} Without loss of generality, let us consider the derivative at time $t=0$ and assume $N\geq 2$ (since 
$N=1$ is a trivial case). First assume that $A(0)=I$, where $I$ is the identity matrix. Then, 
\begin{align*}
\ddd_{t=0} \det(A(t))&=\lim_{t\rightarrow 0}\frac{\det\big[I+ t A'(0)+o(t)\big]-\det(I)}{t}\\
&=\lim_{t\rightarrow 0}\frac{ \Tr\big[I+tA'(0)\big]+o(t)-\Tr(I)}{t}=\Tr(A'(0)).
\end{align*}
In the case where $A(t)$ is invertible at $t=0$, we write 
$$\det(A(t))=\det(A(0)) \det(A(0)^{-1}A(t))$$
and apply the above computation to $\tilde A_h(t)=A(0)^{-1}A(t)$.

For any invertible $A$, we have $\text{\rm com}(A)^t=\det(A)A^{-1}$. Thus, we obtain the last statement by 
extending the formula $\partial_t \det(A(t))=\Tr(\text{\rm com}(A(t))^t \partial_tA(t))$ by a density argument. 
\end{demo}

\subsection{Right-inverse of the divergence}\label{section-inv-div}
Let $k\geq 1$, let $\alpha\in (0,1)$ and let $\Omega_0$ be a $\Cc^{k+1,\alpha}-$domain of $\RR^N$. We define
$$X^{k,\alpha}_0=\{u\in \Cc^{k,\alpha}(\overline\Omega_0,\RR^N)\,,\,u=0\text{ on }\partial\Omega_0\},$$
$$Y^{k-1,\alpha}_m=\Big\{v\in \Cc^{k-1,\alpha}(\overline\Omega_0,\RR^N)\,,\,\int_{\Omega_0} v=0\,\Big\}$$
and
$$L~:~u\in X^{k,\alpha}_0~\longmapsto~ \div u \in Y^{k-1,\alpha}_m~.$$ 
Notice that $L$ is well defined since $\int_{\Omega_0} \div u=0$ if $u_{\partial_{\Omega_0}}\equiv 0$. It is shown 
in \cite[Theorem 30]{Dacorogna-Moser} that the operator $L$ admits a bounded linear right-inverse 
$$L^{-1}:Y^{k-1,\alpha}_m\rightarrow X^{k,\alpha}_0$$
that is $LL^{-1}=\id$ and there exists $K>0$ such that $\|L^{-1}v\|_{\Cc^{k,\alpha}}\leq K 
\|v\|_{\Cc^{k-1,\alpha}}$.

This kind of result is classical, in particular in the framework of Sobolev spaces. In addition to 
\cite{Dacorogna-Moser} see also \cite{BB,BC,Dacorogna}.

\subsection{Fixed point theorem with parameter}
Even though the Banach fixed point theorem is long-established, in this paper we need its extension to 
the case where the contraction depends on a parameter. Of course, this extension is also very classical. We 
briefly recall it for sake of completeness in order to detail the problem of the regularity.
\begin{theorem}\label{th_fixed_point}
Let $U$ be an open subset of a Banach space $X$ and let $V$ be an open subset of a Banach space $\Lambda$. Let 
$F\subset U$ be a closed subset of $X$ and let 
$$\Phi~:~(x,\lambda)\in U\times V~\longmapsto~\Phi(x,\lambda)\in X~.$$
Assume that 
\begin{enumerate}
\item[(i)] For all $\lambda$, $\Phi(\cdot,\lambda)$ maps $F$ into $F$.
\item[(ii)] The maps $\Phi(\cdot,\lambda)$ are uniformly contracting in the sense that there exists $k\in [0,1)$ 
such that 
$$\forall (x,y,\lambda)\in U\times U\times V~,~~\|\Phi(x,\lambda)-\Phi(y,\lambda)\|_X\,\leq\, k \|x -y\|_X~.$$
\end{enumerate}
Then, for all $\lambda\in V$, there exists a unique solution $x(\lambda)$ of $x=\Phi(x,\lambda)$ in $F$. Moreover, 
if $\Phi$ is of class $\Cc^k(U\times V,X)$ with $k\in\NN$, then $x(\lambda)$ is also of class $\Cc^k(V,F)$ (the 
derivatives being understood in the Fr\'echet sense).
\end{theorem}
\begin{demo}
The existence and uniqueness of $x(\lambda)$ correspond of course to the classical Banach fixed point theorem. 
Assume that $\Phi$ is continuous. Then, we write
\begin{align*}
\|x(\lambda)-x(\lambda')\|_X&= \|\Phi(x(\lambda),\lambda)-\Phi(x(\lambda'),\lambda')\|_X\\
&\leq \|\Phi(x(\lambda),\lambda)-\Phi(x(\lambda),\lambda')\|_X + 
\|\Phi(x(\lambda),\lambda')-\Phi(x(\lambda'),\lambda')\|_X\\
&\leq \|\Phi(x(\lambda),\lambda)-\Phi(x(\lambda),\lambda')\|_X + k \|x(\lambda)-x(\lambda')\|_X.
\end{align*}
Since $k<1$ and $\lambda'\mapsto \Phi(x(\lambda),\lambda')$ is continuous, we obtain the continuity of 
$\lambda\mapsto x(\lambda)$. If $\Phi$ is of class $\Cc^k$ with $k\geq 1$, then we apply the implicit function 
theorem to the equation $\Psi(x,\lambda)=0$ with $\Psi(x,\lambda)=x-\Phi(x,\lambda)$. Notice that, due to the 
contraction property, $\|D_x\Phi(x,\lambda)\|_{\Lc(X)}\leq k$ and thus $D_x\Psi$ is invertible everywhere.
\end{demo}

\subsection{The flow of a vector field on a compact domain}\label{sect_regflow}
Let $d\geq 1$, $r\geq 0$ and $p\geq 1$. Let $\overline\Omega_0$ be a compact smooth domain of $\RR^N$ and $I$ be a 
compact interval of times. Let $(t,y)\in I\times\overline\Omega_0\longmapsto U(t,y)\in\RR^N$ a vector field which 
is of class $\Cc^r$ in time and $\Cc^p$ in space, meaning that all derivatives $\partial_{t}^{r'}\partial_y^{p'}U$ 
exist and they are continuous for all $r'\leq r$ and $p'\leq p$. We also assume that $U(t,y)=0$ on 
$\partial\Omega_0$. 

We define $t\mapsto\psi(t,x)$ as the flow corresponding to the ODE
\begin{equation}\label{eq_regflow}
\psi(t_0,x)=x~~\text{ and }~~\partial_t \psi(t,x)=U(t,\psi(t,x))~~~~t\in I~.
\end{equation}
Notice that $\psi$ is locally well defined because $(t,y)\mapsto U(t,y)$ is at least lipschitzian in space and 
continuous in time. Moreover, the trajectories $t\in I\mapsto \psi(t,x)\in\overline\Omega_0$ are globally defined 
because $U(t,y)=0$ on the boundary of $\Omega_0$, providing a barrier of equilibrium points. 

The purpose of this appendix is to show the following regularity result. It is of course a well known property. 
However, the uniform bounds are often not stated explicitly and that is why we quickly recall here the arguments 
to 
obtain them. 
\begin{prop}\label{prop_regflow}
Let $r\geq 0$ and $p\geq 1$. If $U$ is of class $\Cc^r$ in time and $\Cc^{p}$ in space, then the flow $\psi$ 
defined by \eqref{eq_regflow} is of class $\Cc^{r+1}$ in time and $\Cc^{p}$ in space. Moreover, for all $M>0$, 
there exist $C(M)>0$ and $\lambda(M)>0$ such that, if $U$ satisfies 
$$\forall t\in I~,~~\|U(t,\cdot)\|_{\Cc^{p}}~\leq ~ M,$$
then 
$$\forall t\in I~,~~\|\psi(t,\cdot)\|_{\Cc^{p}}~\leq ~ C(M)e^{\lambda(M)|t-t_0|}~.$$
\end{prop}
\begin{demo}
The fact that the $\Cc^0$-bound on $U$ yields a $\Cc^0-$bound on $\psi$ simply comes from \eqref{eq_regflow}. 
If $U(t,y)$ is of class $\Cc^1$ with respect to $y$, it is well know that $\psi(t,y)$ is a of class $\Cc^1$ with 
respect to $y$ and the derivatives solve the ODE
\begin{equation}\label{eq_regflow2}
\partial_t(\partial_{y_i}\psi(t,y))=D_{y}U(t,\psi(t,y)).\partial_{y_i}\psi(t,y)~,
\end{equation}
see for example \cite{Hale}. We have $\partial_{y_i}\psi(t=0)\equiv 0$, thus \eqref{eq_regflow2} and Gr\"onwall's 
inequality ensures the bound on $\partial_{y_i}\psi$.

If $U$ is of class $\Cc^2$ in $y$, the above arguments show that $y\mapsto D_yU(t,\psi(t,y))$ is of class $\Cc^1$ 
and we apply again the procedure to \eqref{eq_regflow2}. We obtain that $\psi(t,y)$ is a of class $\Cc^2$ with 
respect to $y$ and the derivatives solve the ODE
$$\partial_t(\partial^2_{y_iy_j}\psi(t,y))=D_yU(t,\psi(t,y)).\partial^2_{y_iy_j}\psi(t,y)+D^2_{yy}U(t,\psi(t,
y)).(\partial_{y_i}\psi(t,y),\partial_{y_j}\psi(t,y))~,$$
where $\partial^2_{y_iy_j}\psi(t,y)$ is the unknown. Since we already have bounds on $\psi$ and its first 
derivatives, again, Gr\"onwall's inequality yields the bound on $\partial^2_{y_iy_j}\psi$.

By applying the argument as many times as needed, we obtain the uniform bounds for all the wanted derivatives.
We also proceed in the same way to obtain the regularity with respect to the time $t$.
\end{demo}

\subsection{Globalization of local diffeomorphisms}\label{appendix_MO}
In this appendix, we consider a $\Cc^1-$function $\varphi$ for a domain $\Omega_0$ into itself such that $\varphi$ 
is a local diffeomorphism and $\varphi_{|\partial\Omega}=\id$. We would like to obtain that $\varphi$ is in fact a 
global diffeomorphism from $\Omega_0$ into itself. This extension needs topological arguments contained in the 
article \cite{Meisters-Olech} of Meisters and Olech. 

Theorem 1 of \cite{Meisters-Olech} applied to the ball of $\RR^N$ writes as follows.
\begin{theorem}\label{th_Meisters_Olech}
{\bf (Meisters-Olech, 1963)}\\
Let $B_R$ be the open ball of center $0$ and radius $R>0$ of $\RR^N$ and let $\overline B_R$ the closed ball. Let 
$f$ be a continuous mapping of $\overline B_R$ into itself which is locally one-to-one on $\overline B_R\setminus 
Z$, where $Z\cap B_R$ is discrete and $Z$ does not cover the whole boundary $\partial B_R$. If $f$ is one-to-one 
from $\partial B_R$ into itself, then $f$ is an homeomorphism of $\overline B_R$ onto itself. 
\end{theorem}

In fact, the original result of \cite{Meisters-Olech} includes different domains than the balls. Nevertheless, 
if we consider any smooth domain, then it has to be diffeomorphic to a ball (typically, annulus are 
excluded). To consider more general domains, we assume that $f$ is the identity at the boundary.
\begin{theorem}\label{th_MO_2}
Let $\Omega\subset\RR^N$ be a bounded open domain. Let $f$ be a continuous mapping of $\overline \Omega$ into 
itself which is locally one-to-one on $\overline\Omega\setminus Z$, where $Z$ is a finite set. Assume moreover that 
$f$ is the identity on $\partial\Omega$. Then, $f$ is an homeomorphism of $\overline\Omega$ onto itself.   
\end{theorem}
\begin{demo}
For $R$ large enough, $\overline\Omega$ is included inside the ball $B_R$. We extend continuously $f$ to a function 
$\tilde f$ by setting $\tilde f=\id$ on $\overline B_R\setminus\Omega$. Notice that $f$ maps $\overline\Omega$ into 
itself and is locally one-to-one at all the points of the boundary, except maybe at a finite number of them. This 
yields that $\tilde f$ is locally one-to-one at all these points since the extension maps the outside of 
$\overline\Omega$ into itself. We apply Theorem \ref{th_Meisters_Olech} to $\tilde f$ and obtain that $\tilde f$ is 
an homeomorphism of $B_R$. Since it is the identity outside $\Omega$, $f=\tilde f_{|\overline\Omega}$ must be an 
homeomorphism of $\overline\Omega$.  
\end{demo}

If we consider $f$ of class $\Cc^k$ and $Df$ its jacobian matrix, then we may check the local one-to-one property 
by assuming that $\det(Df)$ only vanishes at a finite number of points. More importantly, if $\det(Df)$ never 
vanishes, then $f$ is a diffeomorphism.
\begin{coro}\label{coro_MO}
Let $k\geq 1$, let $\Omega\subset\RR^N$ be a bounded open domain of class $\Cc^k$ and let $f\in\Cc^k(\overline 
\Omega,\overline \Omega)$. Assume that $\det(Df)$ does not vanish on $\overline\Omega$ and that $f$ is the identity 
on $\partial\Omega$. Then, $f$ is a $\Cc^k-$diffeomorphism of $\overline\Omega$ onto itself. 
\end{coro}

We could also be interested in the following other consequence.
\begin{coro}\label{coro_MO_2}
Let $k\geq 1$, let $\Omega\subset\RR^N$ be a bounded open domain of class $\Cc^k$ and let $f$ be a 
$\Cc^k-$diffeomorphism from $\overline \Omega$ onto itself. Assume that $f$ is the identity on $\partial\Omega$.
Then, for all $\varepsilon>0$, there exists $\eta>0$ such that, for all functions $g\in\Cc^k(\overline 
\Omega,\overline \Omega)$ with $g_{\partial\Omega}=\id$ and $\|f-g\|_{\Cc^k}\leq \eta$, $g$ is also a 
$\Cc^k-$diffeomorphism of $\overline\Omega$ onto itself and $\|f^{-1}-g^{-1}\|_{\Cc^k}\leq\varepsilon$.
\end{coro}
\begin{demo}
We simply notice that $\overline\Omega$ is compact and so $|\det(Df)|\geq\alpha>0$ for some uniform positive 
$\alpha$. Thus, for $g$ which is $\Cc^1-$close to $f$, $Dg$ is still invertible everywhere and $g$ is a 
$\Cc^k-$diffeomorphism due to Corollary \ref{coro_MO}.  
\end{demo}

%%%%%%%%%%%%%%%%%%%%%%%%%%%%%%%%%%%%%%%%%%%%%%%%%%%%%%%%%%%%%%%%%%%%%%%%%%%%%%%%
%%%%%%%%%%%%%%%%%%%%%%%%%%%%%%%%%%%%%%%%%%%%%%%%%%%%%%%%%%%%%%%%%%%%%%%%%%%%%%%%
%%%%%%%%%%%%%%%%%%%%%%%%%%%%%%%%%%%%%%%%%%%%%%%%%%%%%%%%%%%%%%%%%%%%%%%%%%%%%%%%
%%%%%%%%%%%%%%%%%%%%%%%%%%%%%%%%%%%%%%%%%%%%%%%%%%%%%%%%%%%%%%%%%%%%%%%%%%%%%%%%

%%%%%%%%%%%%%%%%%%%%%%%%%%%%%%%%%%%%%%%%%%%%%%%%%%%%%%%%%%%%%%%%%%%%%%%
%%%%%%%%%%%%%%%%%%%%%%%%%%%%%%%%%%%%%%%%%%%%%%%%%%%%%%%%%%%%%%%%%%%%%%%

\end{document}